\documentclass[12pt,a4paper]{article}
\usepackage[T2A]{fontenc}
\usepackage[cp1251]{inputenc}
\usepackage[russian, english]{babel}
\usepackage{epsfig,amssymb,color,graphics}
\usepackage{amsthm, amsmath, amssymb, amsbsy, amsfonts,  latexsym, euscript}
\usepackage{tikz}
\usepackage{caption}
\usetikzlibrary{arrows,%
                topaths}%
\usepackage{tkz-berge}
\usepackage{verbatim}
cm \hoffset = -1.0 true cm

\newtheorem{Th}{Theorem}
\newtheorem{st}{Statement}[section]
\newtheorem{Cor}{Corollary}[section]
\newtheorem{Prop}{Proposition}[section]
\newtheorem{Lemm}{Lemma}[section]
\newtheorem{Defin}{Definition}[section]

\newenvironment{demo}{{\bf Proof: }}
{\hfill $\diamond$\medskip}
\begin{document}

\title{
\Large {On Topological Classification of Morse-Smale Diffeomorphisms on the Sphere $S^n$} }
\author{V.~Grines, E.~Gurevich, O.~Pochinka, D.~Malyshev}
\date{National Research University Higher School of Economics,\\  Russian Federation}

\maketitle
\abstract{We consider a class  $G(S^n)$  of  orientation preserving Morse-Smale diffeomorphisms of the sphere  $S^{n}$ of dimension  $n>3$ in assumption that  invariant manifolds of  different saddle periodic points have no intersection.  We put in a correspondence for every diffeomorphism  $f\in G(S^n)$ a colored graph  
$\Gamma_f$ enriched by an automorphism $P_f$. Then we define the notion of isomorphism between two colored graphs and  prove that two diffeomorphisms $f, f'\in G(S^n)$ are  topologically conjugated iff the graphs $\Gamma_f$, $\Gamma_f'$ are isomorphic. Moreover we    establish the existence of a linear-time algorithm for distinguishing two colored graphs of diffeomorphisms from the class $G(S^n)$.}

\sloppy

\section{Introduction and The Statement of Results}

A problem of the topological classification of dynamical systems takes its origin in papers of A.~Andronov, L.~Pontryagin, E.~Leontovich, A.~Mayer and M.~Pexoto. In 1937  A.~Andronov and L.~Pontryagin introduced a notion of the roughness  of  dynamical systems and showed that the necessary and sufficient conditions of the  roughness of a flow on the 2-dimensional sphere are finiteness of non-wandering orbits set, its hyperbolicity and the 
absence of trajectories joining two saddle equilibria or going  from a saddle to the same saddle. In 1960 S.~Smale introduced a similar class of dynamical systems on manifolds of an arbitrary dimension and transform  the condition of the  absence of trajectories between two saddle equilibria into more general condition of transversality of the intersection of invariant manifolds of equilibria and periodic orbits. Later such systems were   called  Morse-Smale systems.
  
   Finiteness of the set of non-wandering orbits leads to an idea of reducing a problem of topological classification of Morse-Smale systems   to a combinatorial task of a description of the  mutual arrangement of such orbits and their invariant manifolds in the ambient  manifolds.   First time this approach was applied by E.~Leontovich and A.~Mayer for classification of flows on two-dimensional sphere with finite set of singular orbits and was developed in papers of M.~Peixoto, A. Oshemkov, V. Sharko,   Y. Umanskii, S. Pilyugin where a similar problem  was solved for Morse-Smale flows on closed manifolds of dimension 2, 3 and higher, and also by V.~Grines and A.~Bezdenezhnych for Morse-Smale diffeomorphisms with finite number of heteroclinic orbits on surfaces\footnote{Non-empty  intersection of invariant manifolds of different saddle periodic points is called  a heteroclinic intersection, an isolated point of such intersection is called a heteroclinic point and its orbit is called a heteroclinic orbit.}. 

It was tuned that this idea, in general, does not work in case of diffeomorphisms on three-dimensional manifolds due to possibility of wild embedding of closures of invariant manifolds of saddle periodic points. This fact required a new language for obtaining of topological invariants. Complete topological classification of Morse-Smale diffeomorphisms on three-dimensional manifolds was obtained in  a series of papers by  Ch.~Bonatti, V.~Grines, F. Laudenbach, O.~Pochinka, E. Pecou,  and V.~Medvedev   (see a reviews~\cite{GrPo} and~\cite{Duke} for references). New invariant called  a scheme of diffeomorphism includes a  topological structures of the   orbit space of action of the diffeomorphism on some wandering set   and embedding   of  projections of invariant manifolds in this orbit space.

Surprisingly it turned out that for some classes of Morse-Smale cascades  on manifolds of  dimension greater than three  a complete invariant  can be borrowed from the flows.  In papers~\cite{GrGuMe08}, \cite{GrGuMe10} by V.Grines, E. Gurevich and V. Medvedev  was obtained  the topological classification of Morse-Smale diffeomorphisms  on manifolds of dimension $n\geq 4$ under  suggestions that all unstable manifolds of saddle periodic points have dimension one and wandering set of diffeomorphisms do not contain heteroclinic orbits. The complete invariant for such  diffeomorphisms is a graph equipped by an automorphism  similar to the Leontovich-Mayer-Peixoto graph in the  two-dimensional case.

The present paper is a continuation of articles~\cite{GrGuMe08}, \cite{GrGuMe10}. We consider the  class  $G(S^n)$  of preserving orientation  Morse-Smale diffeomorphisms on the sphere $S^n$  ($n\geq 4$)  such that the stable and unstable manifolds of different saddle periodic points of any diffeomorphism  from $G(S^n)$ have no intersections. For diffeomorphisms from  the class $G(S^n)$ we provide a combinatorial invariant called the colored graph. For its  exact definition recall some facts.

Let   $\Omega_f$ be a  non-wandering set of a diffeomorphism  $f\in G(S^n)$ and   $\Omega^i_f=\{p\in \Omega_f|\  dim~W^u_p=i\}$, $i\in \{0,1,\dots,n\}$. It follows from~\cite[Theorem~3]{GrGuPo-matan} (see also, \cite[Proposition 4.2]{GrGuPo-emb}) that for $f\in G(S^n)$ the  sets $\Omega^j_f$ are empty for $j\in\{2,\dots,n-2\}$.  

	Let $p \in \Omega^1_f$.  It follows from~\cite[Theorem~2.3]{Sm}  that   the  closure $cl\,W^{s}_p$ of its invariant manifold  $W^{s}_p$ contains, apart the $W^{s}_p$ itself, exactly one periodic point and this point is a source $\alpha$.  Then   the set   $cl~W^{s}_p$ is homeomorphic to  the  sphere of dimension  $(n-1)$ and $cl~W^{s}_p$ is smoothly embedded  at all points except  possibly the  point $\alpha$.  J. Cantrell proved in~\cite{Ca63} that $(n-1)$-dimensional sphere $S^{n-1}\subset S^n$ cannot have one point of wildness (that contrasts with the case $n=3$). Hence the sphere $cl~W^{s}_p$ is locally flat at every point\footnote{A manifold   $N^{k}\subset M^n$ of dimension  $k$ without boundary  is  {\it locally flat in a point $x\in N^k$} if there exists a neighborhood 
$U(x)\subset M^{n}$ of the point    $x$ and a homeomorphism  $\varphi:U(x)\to
\mathbb{R}^{n}$ such that $\varphi(N^{k}\cap U(x))=
\mathbb{R}^{k}$, where  $\mathbb{R}^k=\{(x_1,...,x_n)\in \mathbb{R}^n|\  x_{k+1}=x_{k+2}=...=x_n=0\}$. If the condition of local flatness fails in a point  $x\in N^k$ then the manifold  $N^k$ is called  {\it wild} and  $x$ is called {\it a point of wildness}.}, and, due to~\cite[Theorem 4]{Br62} and \cite[Theorem 5]{Br60} cuts the ambient sphere $S^n$ into two connected components, the closure of each of which is the  ball.  

	
Denote by  $\mathcal{L}_f$ the set of all  the spheres $\{cl~W^s_p, p\in \Omega_f^1\}$ and the spheres $\{cl~W^u_q, q\in \Omega_f^{n-1}\}$  and  put $k_f=|\Omega_f^1\cup \Omega_f^{n-1}|$ (here $|X|$ means cardinality of a  set $X$). Since each $(n-1)$-sphere from the set $\mathcal{L}_f$ cuts the sphere  $S^n$ into two connected components, the set    $S^n\setminus (\bigcup\limits_{p\in \Omega_f^1}cl~W^s_p \cup \bigcup\limits_{q\in \Omega_f^{n-1}}cl~W^u_q)$ consists of  $k_f+1$ connected  components  $D_1, \dots, D_{k_f+1}$. Denote by $\mathcal{D}_f$ the set of all these components.   

  \begin{Defin} A   colored graph of the diffeomorphism  $f\in G(S^n)$ is a graph $\Gamma_f$  with the following properties:
\begin{enumerate}
\item[$1)$]  the set $V(\Gamma_f)$ of vertices of the graph  $\Gamma_f$ is isomorphic to the set  $\mathcal{D}_f$,  the set $E(\Gamma_f)$  of edges of the graph   $\Gamma_f$ is isomorphic to the set  $\mathcal{L}_f$;  
\item[$2)$] vertices $v_i, v_j$ are incident to an  edge  $e_{i, j}$ if and only if correspondent   domains  $D_i, D_j$ have a common boundary; 
\item[$3)$] an edge  $e_{i, j}$ has a color  $s$ ($u$) if it corresponds to a manifold  $cl\,W^s_p\subset \mathcal{L}_f$ ($cl\,W^u_q \subset \mathcal{L}_f$).
\end{enumerate}
\end{Defin}

\begin{figure}
\begin{center}\vspace{1cm}
\includegraphics[width=0.8\linewidth]{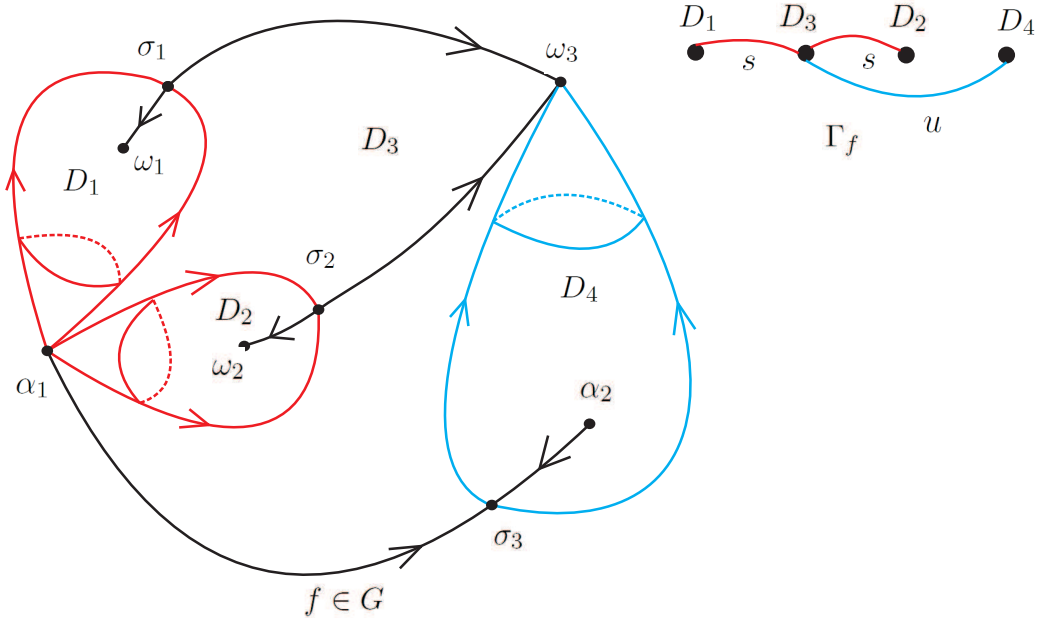}
\caption{A phase portrait of a  diffeomorphism  $f\in G(S^n)$  and its colored graph $\Gamma_f$}\label{graph}
\end{center}
\end{figure}

On the  figure~\ref{graph} is shown a phase portrait of a diffeomorphism $f\in G(S^n)$ with non-wandering set consisting of eight points: sources $\alpha_1, \alpha_2$, saddles $\sigma_1,\sigma_2, \sigma_3$ and sinks $\omega_1,\omega_2,\omega_3$, and its		 colored graph $\Gamma_f$.

Denoty by $\xi: V(\Gamma_f)\to \mathcal{D}_f$ an arbitrary isomorphism and define  an automorphism $P_f: V(\Gamma_f)\to V(\Gamma_f)$ by $P_f=\xi^{-1}f\xi|_{_{V(\Gamma_f)}}$.

\begin{Defin} Graphs $\Gamma_f, \Gamma_{f'}$ of diffeomorphisms $f,f'\in G(S^n)$ are isomorphic if there exists an isomorphisms $\zeta: V(\Gamma_f)\to V(\Gamma_{f'})$  preserving color of edges such that $P_{f'}=\zeta P_f \zeta^{-1}$.
\end{Defin}

\begin{Th}\label{mainth}
Diffeomorphisms   $f,f'\in G(S^n)$ are topologically conjugated if and only if their graphs  $\Gamma_f,  \Gamma_{f'}$  are isomorphic. \end{Th}

The following result shows that the colored graph is the most effective invariant for classification of diffeomorphisms from the class  $G(S^n)$ because there exists an optimal linear algorithm of recognizing such graphs up   isomorphism.

\begin{Th}\label{dima} Let  $\Gamma_f, \Gamma_{f'}$ be graphs of diffeomorphisms  $f, f'\in G(S^n)$  with the same number  $k$ of vertices.  Then there exists an algorithm of a verification of an existence of  an isomorphism between  $\Gamma_f, \Gamma_{f'}$ for the time $O(k)$.
\end{Th}

The structure of the paper is the following. In Section~\ref{start}  we list all auxiliary  results necessary for proof of the Theorem~\ref{mainth}, then provide the proof of Theorem~\ref{mainth} in section~\ref{start1}. Proofs of the auxiliary results required special techniques   are  given in Sections~\ref{emb}-\ref{orbits1}.   The proof of Theorem~\ref{dima} is  given  in  Section~\ref{dima0} that can be read independently from other sections.

\medskip  
{\bf{Acknowledgments}}

Research   is  supported by Russian Science Foundation, project 17-11-01041, except the sections~\ref{start1},~\ref{dima0}  that are supported by Laboratory of Dynamical
Systems and Applications  of National Research University Higher School of Economics, grant   of the Ministry of Science and Higher Education of
the Russian Federation  № 075-15-2019-1931 and by LATNA laboratory of National Research University Higher School of Economics.

\section{Auxiliary results}
\label{start}

\subsection{Linearizing neighborhood and canonical manifolds connected with hyperbolic periodic points}\label{canon-intro}

Define a homeomorphism   $b_{\nu}: \mathbb{R}^{n}\to \mathbb{R}^{n}$, $\nu \in\{+1,-1\}$, by $$b_{\nu}({x}_{1},x_2,\ldots,x_n)=\Big(\nu 2x_{1}, \frac{1}{2}x_{2},\ldots, \frac{1}{2}x_{n-1}, \nu \frac{1}{2} x_{n} \Big).$$ The Origin  $O$ is the unique fixed point of  $b_{\nu}$,  and it is a hyperbolic saddle point. The stable manifold  $W^s_{_{O}}$ coincides with the hyperplane  $x_{1}=0$ and the unstable manifold  $W^u_{_{O}}$ coincides with the $Ox_{1}$-axes.

Put $\mathbb{U}=\{(x_{1},\ldots,x_{n})\in
\mathbb{R}^{n} \mid x_{n}^{2}(x_{1}^{2}+\ldots+x_{n-1}^{2}) \leq 1\}$, $\mathbb{U}_0=\{(x_1,\ldots,x_n)|\ x_1=0\}$, $\mathbb{N}^u=\mathbb{U}\setminus \mathbb{U}_0$, $\mathbb{N}^s=\mathbb{U}\setminus Ox_1$, $\mathbb{K}^{n-1}_{\nu}=\mathbb{U}_0/_{b_\nu}$,
 $\widehat{\mathbb{N}}^s_\nu=\mathbb{N}^s/_{b_\nu}$,
$\widehat{\mathbb{N}}^u_\nu=\mathbb{N}^u/_{b_\nu}$ and denote by $p^s_{b_\nu}:  {\mathbb{N}}^s_\nu\to \widehat{\mathbb{N}}^s_\nu$, $p^u_{b_\nu}:  {\mathbb{N}}^u_\nu\to \widehat{\mathbb{N}}^u_\nu$ the canonical projections. 


Recall that  {\it an  $n$-ball $(n$-disk$)$} is a manifold $B^n$ homeomorphic to the unit ball  $\mathbb{B}^{n}=\{(x_{1},...,x_{n})\in\mathbb R^{n} \mid
x_{1}^{2}+...+x_{n}^{2}\leq1\}$, $n\geq 1$. An open  $n$-ball ($n$-disk) is a manifold homeomorphic to the interior of  $\mathbb{B}^{n}$ and a  sphere  $S^{n-1}$ is a manifold homeomorphic to the boundary  $\mathbb{S}^{n-1}$ of the ball  $\mathbb{B}^n$. The sphere $S^1$ also is called  {\it the knot}.

It is  easy to show that the  manifold  $\mathbb{K}^{n-1}_{+1}$
is homeomorphic to $\mathbb{S}^{n-2}\times \mathbb{S}^1$.  Call the manifold $\mathbb{K}^{n-1}_{-1}$ {\it the standard generalized  $(n-1)$-dimensional Klein Bottle} and a   manifold homeomorphic to 
 $\mathbb{K}^{n-1}_{-1}$  {\it  the generalized  Klein Bottle}. A canonical projection  $p^s_{b_{-1}}|_{\mathbb{U}_0}:\mathbb{U}_0\to \mathbb{K}^{n-1}_{-1}$ induces on  $\mathbb{K}^{n-1}_{-1}$  a structure of non-oriented fiber bundle  over $\mathbb{S}^1$ with a fiber $\mathbb{S}^{n-2}$.  Hence 
 $K^{n-1}_{-1}$ is a non-oriented manifold. Since  $\mathbb{U}_0$ is the universal cover for  $\mathbb{K}^{n-1}_{\nu}$ then the fundamental group  $\pi_{1}(\mathbb{K}^{n-1}_{\nu})$  is isomorphic to $\mathbb{Z}$ (see \cite[Corollary  19.4]{Ko}).

 Call the  manifold $\widehat{\mathbb{N}}^s_{-1}$ {\it the canonical neighborhood of the generalized  Klein Bottle} $\mathbb{K}^{n-1}_{-1}$. 

The proposition below follows directly from the definition of $\widehat{\mathbb{N}}_{\nu}$.

\begin{Prop}
\label{st-nbh}
$ $
\begin{enumerate}
\item   $\widehat{\mathbb{N}}^u_{+1}$ consists of two connected components each of which is diffeomorphic to the direct product  $\mathbb{B}^{n-1}\times \mathbb{S}^1$.

\item   $\widehat{\mathbb{N}}^u_{-1}$  is diffeomorphic  to  $\mathbb{B}^{n-1}\times \mathbb{S}^1$.

\item  $\widehat{\mathbb{N}}^s_{+1}$ is diffeomorphic to   $\mathbb{K}^{n-1}_{+1}\times [-1,1]$.

\item  $\widehat{\mathbb{N}}^s_{-1}$  is a tubular neighborhood of a zero section of non-orientable one-dimensional vector bundle over 
 $\mathbb{K}^{n-1}_{-}$, 
the boundary $\partial \widehat{\mathbb{N}}_{-}$ of  $\widehat{\mathbb{N}}_{-}$ is diffeomorphic to  $\mathbb{S}^{n-2}\times \mathbb{S}^1$, and if 
 $i_{*}: \pi_1(\partial
\widehat{\mathbb{N}}_{-})\to \pi_1(\widehat{\mathbb{N}}_{-})$ is a homomorphism induced by an inclusion then 
$\eta_{_{\widehat{\mathbb{N}}_{-}}}(i_*(\pi_1(\partial \widehat{\mathbb{N}}_{-})))=2\mathbb{Z}$\footnote{Recall that  {\it fiber bundle}  is a structure   $\xi=\{E,B,Y,\pi\}$ where  $E,B,Y$ are topological spaces and  $\pi:E\to B$ is a continuous map such that    for every point $x\in E$ there is an open neighborhood  $U\in B$ of $\pi(x)$ and a homeomorphism  $\varphi:\pi^{-1}(U)\to U\times Y$  such that if  $p_1: U\times Y\to U$ is the projection on the first factor  ($p_1(x,y)=x$) then  $\pi|_{_{\pi^{-1}(U)}}= p_1\varphi|_{_{\pi^{-1}(U)}}$.  $E, B$ and  $Y$ are called  \textit{the total  space, the base  space} and 
\textit{the fiber}; a pair $(U,\varphi)$ is called a chart,  the set 
$\{(U,\varphi)\}$ is called  \textit{a local trivialization} of the bundle. For every loop  $\lambda\subset B$ with the endpoints in a point  $x$ a mapping class group   $T_\lambda:\xi_x\to \xi_x$   induced by chart transition maps along the loop is called {\it a monodromy action}.

\textit{Vector bundle of dimension $n$} is a fiber bundle   $\xi=\{E,B,\mathbb{R}^n,\pi\}$ such that for any two charts  $(U,\varphi),\ (V,\psi)$ such that $x\in U\cap V$ the following condition holds: if  $\varphi_x=p_2\varphi|_{\pi^{-1}(x)}$, $\psi_x=p_2\psi|_{\pi^{-1}(x)}$ then the map  $\psi_x^{-1}\varphi_x: \mathbb{R}^n\to \mathbb{R}^n$  is linear (here  $p_2: U\times \mathbb{R}^n\to \mathbb{R}^n$  is the projection onto the second factor). The fiber $\xi_x=\pi^{-1}(x)$  over the  point  $x\in B$  is enriched by a structure of vector space such that the map  $\psi_x: \xi_x\to \mathbb{R}^n$ is an isomorphism of vector spaces. A zero section of vector bundle is an image  $\zeta(B)\subset E$  of a map  $\zeta:B\to E$ which maps the point $x\in B$  into the  zero vector of the space  $\xi_x$.}.
\end{enumerate}
\end{Prop}

\begin{figure}
\begin{center}\vspace{1cm}
\includegraphics[width=0.8\linewidth]{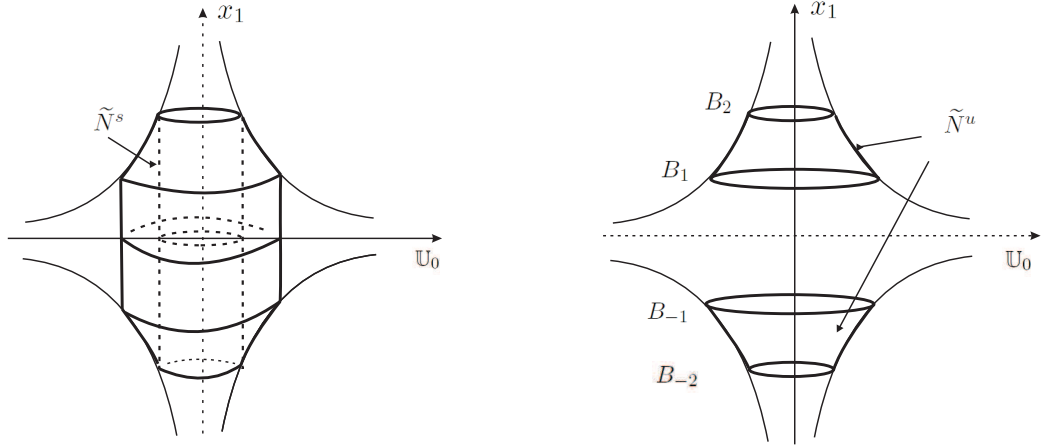}
\caption{Fundamental domains   $\widetilde{N}^s, \widetilde{N}^u$ of action  of the\\ diffeomorphism   $b_{+1}$ on the sets  $\mathbb{N}^s, \mathbb{N}^u$}
\label{canon-nbh}
\end{center}\vspace{1cm}
\end{figure}

On figure~\ref{canon-nbh} the neighborhoods  $\mathbb{N}^s, \mathbb{N}^u$ and fundamental domains   $\widetilde{N}^s=\{(x_1,\dots,x_n)\in \mathbb{N}^s|  \frac14\leq x_2^2+\dots+x_n^2\leq 1\}, \widetilde{N}^u=\{(x_1,\dots,x_n)\in \mathbb{N}^u| |x_1|\in [1,2]\}$ of action of the diffeomorphism   $b_{+1}$ on them are shown\footnote{A fundamental domain of action of a group    $G$ on the set   $X$  is a closed set   $D_{_{G}}\subset X$ containing the subset    $\tilde{D}_{_{G}}$   with the following properties: 1) $cl\,\tilde{D}_{_{G}}=D_{_{G}}$; 2) $g(\tilde{D}_{_{G}})\cap \tilde{D}_{_{G}}=\emptyset$ for any   $g\in G$ different from the neutral element; 3) $\bigcup\limits_{g\in G}g(\tilde{D}_{_{G}})=X$.}. Put   $\mathcal{C}=\{\{(x_1,\dots,x_n)\in \mathbb{R}^n|  \frac14\leq x_2^2+\dots+x_n^2\leq 1\}$.   The set   $\mathbb{N}^s$  is a union of hyperplanes    $\mathcal{L}_t=\{(x_1,...,x_n)\in {N}^s| x_1^2(x_2^2+\dots+x_n^2)=t^2\}, t\in [-1,1]$.  Then the fundamental domain   $\widetilde{{N}}^{s}_{+1}$ is the union of pairs of annuli      $\mathcal{K}_t=\mathcal{L}_t\cap \mathcal{C}, t\in [-1,1]$  and the space     $\widehat{\mathbb{N}}^{s}_{+1}$ can be obtained from    $\widetilde{{N}}^{s}$ by gluing of connected components of the boundary of each annulus by means of  the diffeomorphism    $b_{+1}$.  
The set    $\widetilde{N}^u$  consists of two connected components each of which is homeomorphic to   $\mathbb{B}^{n-1}\times [0,1]$.  The space   $\widehat{\mathbb{N}}^{u}_{+1}$  can be obtained from   $\widetilde{N}^u$  by gluing a disk    $B_1=\{(x_1,\dots,x_n)\in \mathbb{N}^u| x_1=1\}$  with a disk    $B_2=\{(x_1,\dots,x_n)\in \mathbb{N}^u| x_1=2\}$ and a disk    $B_{-1}=\{(x_1,\dots,x_n)\in \mathbb{N}^u| x_1=-1\}$  with a disk   $B_{-2}=\{(x_1,\dots,x_n)\in \mathbb{N}^u| x_1=-2\}$  by means of    $b_{+1}$.

A fundamental domain of action of the diffeomorphism $b_{-1}$ on the set  $\mathbb{N}^u$  is a set   $\widetilde{N}^u_{-1}=\{(x_1,\dots,x_n)\in \mathbb{N}^u| |x_1|\in [1,4]\}$
The space   $\widehat{\mathbb{N}}^{u}_{-1}$ can be obtained from     $\widetilde{N}^u_{-1}$  by gluing a disk  $B_1=\{(x_1,\dots,x_n)\in \mathbb{N}^u| x_1=1\}$ with a disk    $B_4=\{(x_1,\dots,x_n)\in \mathbb{N}^u| x_1=4\}$  by means of  $b_{-1}$. The structure of vector bundle on the space  $\widehat{\mathbb{N}}^s_{-1}$ is defined by the natural projection of the one-dimensional foliation of the  set   $\mathbb{N}^s$ by straight lines parallel to the   $Ox_1$-axis.  This bundle is non-orientable as a wind along a  loop  $p^s_{b_\nu}(l_{\nu})$ where $l_{\nu}$ is a segment of  $0x_1$ with endpoints  $(0,\ldots,0,1)$,  $(0,\ldots,
1/2)$ induces  reversing orientation map.


Now  came back to diffeomorphisms from the class $G(S^n)$. For  saddle point  $\sigma\in  \Omega^1_f\cup \Omega^{n-1}_f$ denote by   $m_\sigma$  its period. 

Say that the saddle point $\sigma$ of period $m_\sigma$ has {\it an orientation type}  $\nu_{\sigma}=+1$ ($\nu_\sigma=-1$) if the restriction $f^{m_\sigma}|_{W^u_\sigma}$  preserves (reverses) an orientation of $W^u_\sigma$. 

Due to~\cite[Theorem 2.1.2]{GrPo-book} (see also ~\cite[Proposition 4.3]{GrGuPo-emb}) the following proposition holds. 

\begin{Prop} 
\label{grgupo-ado}
For any diffeomorphism   $f\in G(S^n)$  there exists a set of pairwise  disjoint neighborhoods    $\{N_\sigma, \sigma\in \Omega^1_f\cup \Omega^{n-1}_f\}$ such that for any  $N_\sigma$  there exists a homeomorphism  $\chi_{\sigma}:N_{\sigma} \to \mathbb{U}$ such that 

{\rm 1)}~if $\sigma\in \Omega^1_f$
then  $\chi_{\sigma}f^{m_{\sigma}}|_{_{N_{\sigma}}}=b_{\nu_{\sigma}}\chi_{\sigma}|_{_{N_{\sigma}}}$,

{\rm 2)}~if $\sigma\in \Omega^{n-1}_f$ then
$\chi_{\sigma}f^{m_{\sigma}}|_{_{N_{\sigma}}}=b^{-1}_{\nu_{\sigma}}\chi_{\sigma}|_{_{N_{\sigma}}}$.
\end{Prop}

Call the neighborhood  $N_\sigma$ {\it a  linearizing neighborhood}.

\subsection{The scheme $S_f$ of the diffeomorphism $f\in G(S^n)$}

Here we introduce a notion of  {\it the scheme} of a diffeomorphism $f\in G(S^n)$ which is an effective tool for studying of the dynamics.  Moreover, due to~\cite[Theorem 1]{GrGuPo-matan} the scheme is complete invariant in class $G(S^n)$, so the problem of topological classification is reduced to the proof of the fact that   diffeomorphisms having isomorphic colored graphs  have equivalent schemes too. This fact will be  given in    Lemma~\ref{gr=sh}. In section~\ref{mat} we adduce  an adaptation of the proof of~\cite[Theorem 1]{GrGuPo-matan} for the  class $G(S^n)$ to provide  possibility of independent reading of this paper. 

 
Represent the sphere   $S^n$  as the union of pairwise disjoint sets  $$A_f=(\bigcup\limits_{\sigma\in \Omega^{1}_f}{W^u_\sigma})\cup \Omega^{0}_f, R_f=(\bigcup\limits_{\sigma\in \Omega^{n-1}_f}{W^s_\sigma})\cup  \Omega^{n}_f, V_f=S^n\setminus(A_f\cup R_f).$$  Similar to~\cite{GrPoZh} one can prove that the sets  $A_f, R_f, V_f$  are connected, the set   $A_f$  is an attractor,  the set  $R_f$ is a repeller\footnote{A set  $A$ is called an attractor of a homeomorphism  $f:M^n\to M^n$ if there exists a closed neighborhood  $U\subset M^n$ of the set  $A$ such that $f(U)\subset int~U$ and  $A=\bigcap\limits_{n\geq 0} f(U)$.  A set  $R$ is called a repeller of a homeomorphism  $f$ if it is an attractor for the homeomorphism $f^{-1}$.} and $V_f$ consists of wandering orbits of $f$ moving from  $R_f$ to $A_f$. 

Denote by  $$\widehat V_f=V_f/f$$  the orbit space of the action of  $f$ on $V_f$,  by   $$p_{_f}:V_f\to \widehat V_f$$ the natural projection.  In virtue of~\cite{Th01} (Theorem~3.5.7  and Proposition~3.6.7) $p_{_f}$    is a covering map and the space  $\widehat V_f$ is a manifold.

The following lemma is proved in Section~\ref{orbits1}.

\begin{Lemm}\label{svo1} Let $f\in G(S^n)$. Then 
 $\widehat{V}_f$ is homeomorphic to the direct product $\mathbb{S}^{n-1}\times \mathbb{S}^1$.
\end{Lemm}

Remark that for  $n=3$  Lemma~\ref{svo1} in general   is not true (see, for example, \cite[Section 5]{GrPo-2013}).

Denote by  $$\eta_{_{f}}:\pi_{1}(\widehat V_f)\to \mathbb{Z}$$ a homomorphism defined in the following way. Let $\hat{c}\subset \widehat V_f$ be a loop non-homotopic to zero in   $\widehat V_f$ and   $[\hat{c}]\in \pi_1(\widehat V_f)$ be a homotopy class of  $\hat{c}$.  Choose an arbitrary point  $\hat{x}\in \hat{c}$, denote by   $p_{_{f}}^{-1}(\hat{x})$  the complete inverse image of  $\hat{x}$, and fix a point  $\tilde{x}\in p_{{_{f}}}^{-1}(\hat{x})$. As   $p_{_{f}}$ is the covering map then   there is a unique path  $\tilde {c} (t)$  beginning at the point $\tilde {x}$ ($\tilde
{c} (0) = \tilde x$) and covering the loop $c$ (such that  $p_{_{f}}(\tilde {c} (t))=
\hat{c}$). Then there  exists the element  $n\in\mathbb {Z}$ such that 
$\tilde{c}(1)=f^n (\tilde {x})$. Put  $\eta_{_{f}}([\hat{c}]) =
n$. It follows from~\cite{Ko} (Chapter~18)  that the homomorphism  $\eta_{_{f}}$ is an epimorphism.

Put  

 $$\hat{L}^s_{f}=\{p_f(W^s_p\setminus p),  p\in \Omega^1_f\}),\,\,\, \hat{L}^u_{f}= \{p_{_f}(W^u_q\setminus q), q\in \Omega^{n-1}_f\}. $$ 

For a periodic  saddle point $\sigma \in \Omega^1_f$ ($\sigma \in \Omega^{n-1}_f$)   put $l_\sigma=W^s_\sigma\setminus \sigma$ ($l_\sigma=W^u_\sigma\setminus \sigma$),  $\hat{l}_\sigma=p_f(l_\sigma)$, $\widehat{N}_\sigma=p_f(N_\sigma\cap V_f)$, so $\hat{l}_\sigma$ is an element of the set $\hat{L}^s_{f}\cup \hat{L}^u_{f}$.

Due to~\cite[Theorem~5.5]{Ko}  (see also~\cite[Proposition 1.2.3]{BoGrPo-04}) the superposition $\hat{\chi}_\sigma=p_f\chi_\sigma^{-1}(p^s_{b_{\nu}})^{-1}$ define a homeomorphism from $\mathbb{N}^s_{\nu_\sigma}$ to $\widehat{N}_{\sigma}$  such that $\hat{\chi}_\sigma({\mathbb{K}^{n-1}_{\nu_\sigma}})=\hat{l}_\sigma$    and  the   following corollary from the Proposition~\ref{grgupo-ado} holds. Denote by  $\hat{\chi}_{\sigma*}:\pi_1(\mathbb{K}^{n-1}_{\nu_\sigma})\to \pi_1(\hat{V}_f)$  the homomorphism induced by $\hat{\chi}_\sigma$. 

\begin{Cor}\label{ado-down}    $\eta_{_f} \hat{\chi}_{\sigma*}(\pi_1(\mathbb{K}^{n-1}_{\nu_\sigma}))~=~m_\sigma\mathbb{Z}$.    
\end{Cor}

\begin{Defin}
The collection  $$S_{f}=(\widehat V_{f},\hat{L}^s_{f},\hat{L}^u_{f}, \eta_{_f})$$  is called  {\it the  scheme} of the diffeomorphism   $f\in G(S^n)$.
\end{Defin}

\begin{Defin} Schemes $S_f$ and  $S_{f'}$ of diffeomorphisms $f,f'\in G(S^n)$ are called {\it equivalent} if there exists a homeomorphism  $\hat\varphi:\widehat V_f\to\widehat V_{f'}$ such that   $\hat\varphi(\hat{L}^s_{f})=\hat{L}^s_{f'}$,   $\hat\varphi(\hat{L}^u_{f})=\hat{L}^u_{f'}$, and  $\eta_{_f}=\eta_{_{f'}}\hat\varphi_*$.
\end{Defin}


\subsection{Interrelation between the colored graph and the scheme}
\label{start1}
Here we recall some notion of graph theory, establish  properties of   the colored graphs of diffeomorphisms from class $G(S^n)$ and their  interrelation with  schemes that play key role in the proof of Theorem~\ref{mainth}.

Recall that a loop-less graph without multiple and oriented edges is
said to be \emph{simple}. A \emph{tree} is a connected simple cycle-free graph.  Call a vertex incident to exactly one edge {\it a hanging vertex or a leaf}.  

\begin{Prop}\label{doestree} The colored graph $\Gamma_f$ of a diffeomorphism $f\in G(S^n)$ is a   tree.
\end{Prop}
\begin{demo} By definition  any edge $e$ of graph $\Gamma_f$ corresponds to an $(n-1)$-dimensional sphere which cuts the ambient sphere $S^n$ into two connected components. Then the edge $e$  cuts the  graph $\Gamma_f$ into two connected components, so the graph $\Gamma_f$ does not contains cycles. Prove that the graph is connected.  It is well-known that a connected graph  with $k+1$ vertices  is a tree if and only if  it has $k$ edges (we give the proof of this fact in Proposition~\ref{tr}). If the graph $\Gamma_f$ is disconnected, then it consists of  connected subgraphs $\Gamma_1,\dots, \Gamma_m$, $m\geq 2$. Then adding $(m-1)$ edges makes the union of $\{\Gamma_i\}$ a  connected graph without cycles (that is tree) with $k_f+1$ vertices and $k_f+m$ edges that contradicts to Proposition~\ref{tr}. So, the graph $\Gamma_f$ is connected and does not contains cycles, hence    it is a tree.  
\end{demo}


 Recall that we denoted by $V(\Gamma_f), E(\Gamma)$ the sets of vertices and edges of the graph $\Gamma_f$ correspondingly. Denote by $uv\subset E(\Gamma_f)$ an edge connecting vertices $v,u\in V(\Gamma_f)$. 

We associate with the graph   $\Gamma_f$ a sequence~$\Gamma_{f,0},\Gamma_{f,1},$ $\ldots,\Gamma_{f,s}$ of trees, such that $\Gamma_0=\Gamma_f$, $\Gamma_{f,s}$ contains one or two vertices and, for any $i\in \{1,\dots, s\}$, a tree $\Gamma_{f,i}$ is obtained
from $\Gamma_{f,i-1}$ by deletion of all its leaves. All the vertices of $\Gamma_{f,s}$ are called \emph{central} vertices of the tree $\Gamma_f$ and if $\Gamma_{f,s}$ has an edge then it is called {\it the central} edge of the tree $\Gamma_f$. The tree $\Gamma_f$ will be \emph{central}, if it has exactly one central vertex, and \emph{bicentral}, otherwise. 
\emph{The rank of a vertex $x\in V(\Gamma_f)$}, denoted by $rank(x)$, is the number $\max\{i|~x\in V(\Gamma_{f,i})\}$. 
It follows from the definition that if vertices $v,w$ are incident to a non-central edge than $|rank(v)-rank(w)|=1$, and central vertices of a bicentral tree has the same rank.

 An \emph{automorphism} $P$ of a tree $\Gamma$ is a bijective mapping of~$V(\Gamma)$ onto itself that keeps
adjacency, i.e. $$\forall u,v\in V(\Gamma)~[uv\in E(\Gamma)\Leftrightarrow P(u)P(v)\in E(\Gamma)]. $$ 

An  automorphism $P_f$ can be represented as a superposition of cyclic sub-permutations and the set $V(\Gamma_f)$ can be decomposed into subsets invariant under the  sub-permutations  which are said to be
\emph{orbits}. Clearly that every orbit of~$P_f$ consists of the same rank vertices of $\Gamma_f$ and if the tree is central (bicentral) then its central vertex (central edge) stay fixed under any automorphism.

 Orbits $O_1$ and $O_2$ of $P$ will be called \emph{neighbour}, if there are adjacent vertices $x_1\in O_1$ and $y_1\in O_2$, such that $rank(x_1)-rank(y_1)=1$.\\

\begin{Prop}\label{dima-l1} Let $v_1,w_1$ be  vertices of graph $\Gamma_f$ incident to the edge $v_1w_1$, $rank(v_1)-rank(w_1)=1$ and  $O_1=(v_1,v_2,\ldots,v_p)$ and $O_2=(w_1,w_2,\ldots,w_q)$ be  orbits of $v_1,w_1$ correspondingly such that  $P_f(v_i)=v_{i+1}$, $P_f(w_j)=w_{j+1}$, taking the indices modulo $p$ and $q$, correspondingly. Then the following properties are true:

\begin{enumerate}
\item $q\vdots p$;

\item for any $i\in \overline{1,p}$ the set of neighbours of $v_i$ belonging to $O_2$ coincides with $\{w_i,w_{i+p},w_{i+2p},\ldots,w_{i+(\frac{q}{p}-1)\cdot p}\}$;

\item all edges simultaneously incident to a vertex $O_1$ and a vertex of $O_2$ have the same color.
\end{enumerate}
\end{Prop}
\begin{demo}  Assume that $p>q$. Since $v_1w_1\in E(\Gamma_f)$ and $P_f$ is an automorphism
of $\Gamma_f$  then  for any $i=\overline{1,q}$, we have $v_iw_i\in E(\Gamma_f)$. Therefore $P_f(v_q)P_f(w_q)=v_{q+1}w_1\in E(\Gamma_f)$. Hence the tree $\Gamma_f$ contains a cycle and we have a contradiction. Now, suppose that~$q=l\cdot p+r$, where $l$ is a natural number and
$r$ is a natural number between 1 and $p-1$. Since $v_1w_1\in E(\Gamma_f)$ and $P_f$ is the automorphism  then for any $i=\overline{1,p}$, we have $v_iw_i\in E(\Gamma_f)$.
Therefore $v_rw_{p+r}\in E(\Gamma_f)$. Thus $v_rw_{2p+r},v_rw_{3p+r},\ldots,v_rw_{l\cdot p+r}\in E(\Gamma_f)$. Hence~\mbox{$v_rw_q\in E(\Gamma_f)$}, $v_{r+1}w_1\in E(\Gamma_f)$ and $\Gamma_f$ contains a cycle.  Therefore $q\vdots p$. From this fact and the fact that $P_f$ is the
automorphism of $\Gamma_f$  it follows that  for any $i\in \overline{1,p}$, the set of all neighbours of~$v_i$, belonging to $O_2$, coincides with $\{w_i,w_{i+p},w_{i+2p},\ldots,w_{i+q-p}\}$. As $P_f(v_1)P_f(w_1)=v_2w_2$  then edges $v_1w_1, v_2w_2$ have the same color and it is also true   for  all edges simultaneously incident to a vertex $O_1$ and to a vertex of~$O_2$.  \end{demo}

The automorphism $P_f$ naturally  induces the map of the set $E(\Gamma_f)$ which we will denote by $P_f$ too. Proposition~\ref{dima-l1} immediately leads to the next corollary.

\begin{Cor}\label{dima-c1}
Let $vw\in E(\Gamma_f)$  be a non-central edge and $rank(v)<rank(w)$. Then the  period of the edge $vw$ equals to the period of the vertex $v$. \end{Cor}

Two trees $\Gamma$ and $\Gamma'$ are said to be
\emph{isomorphic}, if there is a bijection  $\xi: V(\Gamma)\longrightarrow V(\Gamma')$, called an \emph{isomorphism}, such that
$\forall u,v\in V(\Gamma)~[uv\in E(\Gamma)\Leftrightarrow \xi(u)\xi(v)\in E(\Gamma')].$

Obviously, under any isomorphism of $\Gamma$, any its vertex must be mapped into a same rank vertex. 

For any $f\in G(S^n)$ define {\it a weighted graph $\widehat{\Gamma}_f$} by  the following: 1) glue vertices in the graph $\Gamma_f$  belonging to the same orbit of the automorphism $P_f$ and correspondent edges; 2)  enrich each new edge of the graph $\widehat{\Gamma}_f$ with the weight equal to the period of the correspondent separatrix of the diffeomorphism $f$.  It follows from Proposition~\ref{dima-l1}  that this gluing operation keeps colors of edges.  If graph $\Gamma_f$ is bicentral and central vertices generates the period 2 orbit then graph $\widehat{\Gamma}_f$ has the  unique  loop corresponding to the central edge of the graph $\Gamma$. Otherwise the graph $\widehat{\Gamma}_f$ is the tree with the same central vertices as $\Gamma_f$.
In both cases we will say that a vertex $\hat{v}\in \widehat{\Gamma}_f$ has rank $k$ if the  rank of the correspondent vertex $v\in \Gamma_f$ equals $k$.

Observation above immediately leads to the  following    proposition. 

\begin{figure}
\begin{center}\vspace{1cm}
\includegraphics[width=0.7\linewidth]{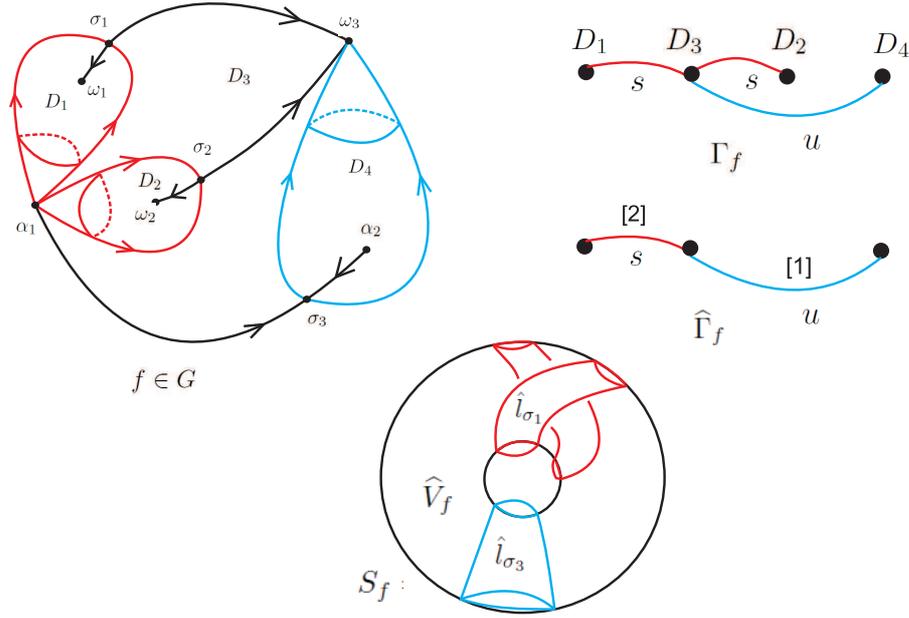}
\caption{A phase portrait of a  diffeomorphism  $f\in G(S^n)$,   its colored graph $\Gamma_f$, weighted graph $\widehat{\Gamma}_f$ and elements of the scheme $S_f$.}\label{graph}
\end{center}
\end{figure}

\begin{Prop}\label{weights} If graphs $\Gamma_f, \Gamma_{f'}$ of diffeomorphisms $f,f'\in G(S^n)$  are isomorphic, than the weighted graphs $\widehat{\Gamma}_f, \widehat{\Gamma}_{f'}$  are also isomorphic  by means of isomorphism preserving weights of edges.
\end{Prop} 

 Denote by $\widehat{\mathcal{D}}_f$ the set of connected components of the set $\widehat{V}_f\setminus (\hat{L}^u_f\cup \hat{L}^s_f)$. Since deleting  $0$-and $1$-dimensional set from manifold of dimension three and greater  does not change the number of connected components, there is natural  one-to-one correspondence $\xi_*$ between the sets $E(\Gamma_f), V(\Gamma_f)$ and $\mathcal{L}_f\cap V_f$, $\mathcal{D}_f\cap V_f$ such that $\xi_*P_f=f\xi_*$.  This fact immediately  provides the following proposition.

\begin{Prop} \label{graph-and-scheme} There is one-to-one correspondences  $\hat{\xi}: E(\widehat{\Gamma}_f)\cup V(\widehat{\Gamma}_f) \to  (\hat{L}^u_f\cup \hat{L}^s_f) \cup \widehat{\mathcal{D}}_f$ such that $\hat{\xi}(E(\widehat{\Gamma}_f))=\hat{L}^u_f\cup \hat{L}^s_f$,  $\hat{\xi}(V(\widehat{\Gamma}_f))= \widehat{\mathcal{D}}_f$.
\end{Prop}

\begin{Prop}\label{svo2} The  set $\Omega^1_f\cup \Omega^{n-1}_f$ contains at most one saddle point of negative orientation type. If such a point $\sigma$ exists then  its period equals 1,  the manifold   $\hat{l}_\sigma$ is homeomorphic to the generalized Klein Bottle and the  $(n-1)$-dimensional invariant  manifold of $\sigma$ corresponds to the  loop of the graph $\widehat{\Gamma}_f$.
\end{Prop}
\begin{demo}

 Suppose that the set  $\Omega^1_f\cup \Omega^{n-1}_f$ contains a point  $\sigma$ of the negative orientation type and period    $m_\sigma$.   For certainty suppose that   $\sigma\subset  \Omega^1_f$. Due to Proposition~\ref{grgupo-ado} there exists the neighborhood   $N_\sigma$ and the homeomorphism  $\chi_\sigma: N_\sigma\to \mathbb{U}$ such that  $f^{m_\sigma}|_{_{N_\sigma}}=\chi_\sigma^{-1}b_{-}\chi_\sigma|_{_{N_\sigma}}$ then one-dimensional separatrices of the point  $\sigma$ have the period  $2m_\sigma$. The closure of the stable manifold of the point  $\sigma$ is locally flat   $(n-1)$-sphere cutting the ambient sphere $S^n$ on two open  $n$-ball  $B_1, B_2$ each of which contains a one-dimensional separatrix of   $\sigma$.  Hence  $f^{m_\sigma}(B_1)=B_2$, $f^{m_\sigma}(B_2)=B_1$.   Hence if we delete from the graph $\Gamma_f$  the  edge  $e_\sigma$ corresponding to the manifold $W^s_{\sigma}$ we get two  graphs isomorphic by means of $P^{m_{\sigma}}_f$, so the edge $e_\sigma$ is the central. Since any tree has no greater than one central edge, there are no more than one saddle point of negative orientation type.   
\end{demo}

For $\widehat{\Gamma}_f$ construct a sequence $\widehat{\Gamma}_{f,0}=\widehat{\Gamma}_f, \widehat{\Gamma}_{f,1},\dots,\widehat{\Gamma}_{f,s}$  similar to the sequence $\Gamma_{f,0},\Gamma_{f,1},\dots,\Gamma_{f,s}$. 

For $i\in \{0,\dots,s-1\}$ put  $\Delta_{f,i}=\hat{\xi}(\widehat{\Gamma}_f\setminus \widehat{\Gamma}_{f,i+1})$. If the graph $\widehat{\Gamma}_f$ have the central edge $e_*$  choose an arbitrary connected component $\widehat{\Gamma}_{f,*}$ of the graph $\widehat{\Gamma}_f$ and put  $\Delta_{f,s}=\hat{\xi}(\widehat{\Gamma}_{f,*}\cup e_*)$. If the graph $\widehat{\Gamma}_f$ have the loop, then there exist a saddle point $\sigma_*\in \Omega^1_f\cup \Omega^{n-1}_f$. Put $\Delta_{f,s}=cl\,(\widehat{V}_f\setminus int\,\hat{N}_{\sigma_*})$, where $N_{\sigma_*}$ is the linearizing neighborhood of the point $\sigma_*$.

It follows from definition that $\Delta_0\subset \Delta_1\subset \dots \subset \Delta_s$.

  The following Lemma is proved in Section~\ref{orbits1}.

\begin{Lemm}\label{svo} Let $f\in G(S^n)$. Then 
 every connected component of the set $\Delta_i$ is homeomorphic to $\mathbb{B}^{n-1}\times \mathbb{S}^1$.  
\end{Lemm}

Now  to  prove that isomorphism of the colored graphs $\Gamma_f, \Gamma_{f'}$ leads to equivalence of the schemes $S_f, S_{f'}$ we need one more  auxiliary proposition. 

Let $M^n$ be a topological manifold, possibly, with a non-empty  boundary $\partial\, M^n$,  $\beta\in M^n$ be a knot, and  $g:\mathbb{B}^{n-1}\times \mathbb{S}^1\to {M}^n$ be a topological embedding such that  $g(\{O\}\times \mathbb{S}^1)=\beta$. The  image $N_\beta=g(\mathbb{B}^{n-1}\times \mathbb{S}^1)$ is called {\it the tubular neighborhood} of the knot $\beta$. 

 The following Proposition is a slight modification of~\cite[Proposition 5.1, 5.3]{GrGuPo-emb}. We prove it in section~\ref{emb}. 

\begin{Prop}\label{nb-equ} Suppose that $n\geq 4$,  $\mathbb{P}^{n-1}$ is either the ball  $\mathbb{B}^{n-1}$ or the sphere  $\mathbb{S}^{n-1}$. Let   $\{\beta_i\}$, $\{\beta'_i\} \subset int\,\mathbb{P}^{n-1}\times \mathbb{S}^1$
are two families of pairwise disjoint  knots such that knots 
$\beta_i,\beta'_i$ are homotopic, $i\in \{1,\dots,k\}$ and  $\{{N}_{\beta_i}\},\{N_{{\beta'_i}}\}\subset \mathbb{P}^{n-1}\times \mathbb{S}^1$ are   pairwise disjoint tubular neighborhoods of knots $\{\beta_i\}, \{\beta'_j\}$.    Then there exists a homeomorphism  $h:\mathbb{P}^{n-1}\times \mathbb{S}^1\to \mathbb{P}^{n-1}\times \mathbb{S}^1$  such that $h(\beta_i)=\beta'_i, h({N}_{\beta_i})={N}_{{\beta'_i}}, i\in \{1,...,k\}$, and if $\mathbb{P}^n=\mathbb{B}^n$ then   $h|_{_{\partial\,\mathbb{P}^{n-1}\times \mathbb{S}^1}}=id$.
\end{Prop}
\begin{Cor}\label{contin} In assumption of Proposition~\ref{nb-equ} let $\tilde{h}:\partial{\mathbb{B}^{n-1}\times \mathbb{S}^1}\to \partial{\mathbb{B}^{n-1}\times \mathbb{S}^1}$ be a homeomorphism. Then there  exists  a homeomorphism  $h:\mathbb{B}^{n-1}\times \mathbb{S}^1\to \mathbb{B}^{n-1}\times \mathbb{S}^1$  such that $h(\beta_i)=\beta'_i, h({N}_{\beta_i})={N}_{{\beta'_i}}, i\in \{1,...,k\}$ and  $h|_{_{\partial\,\mathbb{B}^{n-1}\times \mathbb{S}^1}}=\tilde{h}$. 
\end{Cor}

 \begin{Lemm}\label{gr=sh} If the colored graphs $\Gamma_f, \Gamma_{f'}$ of  diffeomorphisms $f,f'\in G(S^n)$ are isomorphic then the  schemes $S_f, S_{f'}$ are equivalent.  
\end{Lemm}
\begin{demo} Isomorphism $\zeta$ of the graphs  $\Gamma_f, \Gamma_{f'}$  induce isomorphism $\hat{\zeta}: \widehat{\Gamma}_f\to \widehat{\Gamma}_{f'}$ that preserves   weights and ranks  of edges.  Then for any $i\in \{0,\dots,s\}$ isomorphism $\hat{\zeta}$ induce one-to one correspondence $\zeta_i$ between connected components of the sets $\Delta_{f,i}$ and $\Delta_{f',i}$. If graph $\widehat{\Gamma}_f$ is bicentral then  there are two ways to choose the set $\Delta_{f,s}$, but, since graph $\widehat{\Gamma}_f, \widehat{\Gamma}_{f'}$ are isomorphic,   it is possible to choose sets $\Delta_{f,s}, \Delta_{f',s}$ in such a way that for any connected components  $P_{f,i}\in \Delta_{f,i}, P_{f,j}\in \Delta_{f,j}$ such that  $P_{f,i}\subset int\, P_{f,j}$ holds  $\zeta_i(P_{f,i})\subset int\,\zeta_j(P_{f,j})$, $i\in \{0,\dots,s-1\}, j\in \{1,\dots,s\}, i<j$.  

Two cases are possible: 1) there are no saddle points of negative orientation type; 2) there is a saddle point of negative orientation type. 

Consider case 1). Without loss of generality assume that graph $\widehat{\Gamma}_f$ is bicentral. For central graph  arguments are similar.
 
 It  follows from Lemm~\ref{svo} and Proposition~\ref{nb-equ} that   there exists a  homeomorphism  $\psi_s: \hat{V}_f\to \hat{V}_{f'}$ such that $\psi_s(\Delta_{f,s})=\Delta_{f',s}$. If $s=0$ then the proof is complete and $\hat\varphi=\psi_s$.  
		
		 Let $s>0$. Denote the images of  the sets  $\Delta_{f,0},\dots,\Delta_{f, s-1}$  under the homeomorphism  $\psi_s$  by the same symbols as their originals. 
		 
		Due to Corollary~\ref{contin} from Proposition~\ref{nb-equ} there exists a homeomorphism $\psi_{s-1}:\widehat{V}_{f'}\to \widehat{V}_{f'}$ such that  $\psi_{s-1}|_{\widehat{V}_{f'}\setminus int\, \Delta_{f',s}}=id$, $\psi_{s-1}(P_{f,s-1})=\zeta_i(P_{f,s-1})$ for any connected component $P_{f,s-1}\in \Delta_{f,s-1}$.    If    $s=1$  then the prove is complete and  $\hat\varphi=\psi_{s-1}\psi_s$.  Otherwise, continue the process and after finite number of steps get the desired homeomorphism   $\hat\varphi$.
		
 In case 2) denote by $\sigma_*\in \Omega_f$, $\sigma'_*\in\Omega_{f'}$ the points of negative orientation type. Due to Corollary~\ref{ado-down}  there exists a homeomorphism $\psi_*: \hat{N}_{\sigma_*}\to \hat{N}_{\sigma'_*}$. Due to Lemma~\ref{svo} the manifold $cl(\hat{V}_f\setminus \hat{N}_{\sigma_*})$ is homeomorphic to $\mathbb{B}^{n-1}\times\mathbb{S}^1$. Then using Corollary~\ref{contin} similarly to case 1)  it is possible to continue the homeomorphism $\psi_*$ up to the desire homeomorphism of  $\widehat{V}_f$.

\end{demo}     

\section{Topological classification of diffeomorphisms from $G(S^n)$} 
\label{mat}

It is clear that if diffeomorphisms $f,f'\in G(S^n)$ are topologically conjugated then their colored graphs are isomorphic. Prove the inverse.   Let colored graphs $\Gamma_f, \Gamma_f'$ of diffeomorphisms $f,f'\in G(S^n)$ are isomorphic. According to Lemma~\ref{gr=sh} the  schemes $S_f, S_{f'}$ of diffeomorphisms are equivalent, that is  there exists a homeomorphism $\hat\varphi:\widehat V_f\to\widehat V_{f'}$ such that   $\hat\varphi(\hat{L}^s_{f})=\hat{L}^s_{f'}$,   $\hat\varphi(\hat{L}^u_{f})=\hat{L}^u_{f'}$, and  $\eta_{_f}=\eta_{_{f'}}\hat\varphi_*$. Prove that  the diffeomorphisms $f,f'$ are topologically conjugated.

A homeomorphism $\hat\varphi:\widehat V_f\to\widehat V_{f'}$ induced a homeomorphism $\varphi:V_f\to V_{f'}$ such that  $f|_{_{V_f}}=\varphi^{-1}f'\varphi|_{_{V_{f'}}}$  and for any saddle point  $\sigma\in \Omega^1_f$ $(\sigma\in \Omega^{n-1}_f)$  there is a point  $\sigma'\in \Omega^1_{f'}$ $(\sigma'\in \Omega^{n-1}_{f'})$ such that  $\varphi(W^s_{\sigma}\setminus{\sigma})=W^{s}_{\sigma'}\setminus{\sigma'}$ ($\varphi(W^u_{\sigma}\setminus{\sigma})=W^{u}_{\sigma'}\setminus{\sigma'}$). The homeomorphism $\varphi$  extends univocally to the set $\Omega_f$.

Extend the homeomorphism $\varphi$ to the sets $A_f, R_f$.
First, construct a conjugating homeomorphisms $H_1, H_{n-1}$  on the  linearizing neighborhoods $\{N_\sigma\}$ of saddle points from sets $\Omega^1_f, \Omega^{n-1}_f$ that  coincides with $\varphi$ on the boundaries of the linearizing neighborhoods. Then the desire homeomorphism $H:S^n\to S^n$ conjugating $f$ and $f'$ will be defined by 

$$
H(x)=
\begin{cases}
\varphi(x),&
\displaystyle{x\in S^n\setminus
\bigcup\limits_{\sigma \in \Omega_f^1\cup \Omega_f^{n-1}}N_{\sigma},}
\\[5pt]
H_\delta(x),
&
x\in N_\sigma, \sigma\in \Omega^{\delta}_f, \delta\in \{1,n-1\}.
\end{cases}
$$   
 
 To built the homeomorphism $H_1$ choose
 exactly one point of each saddle orbit of index 1 and denote
the obtained  set by $\widetilde{\Omega}^1_f$. Due to proposition~
\ref{grgupo-ado}  there exists a family of pairwise disjointed neighborhoods  $\{N_{\sigma}\} (\{N_\sigma'\})$ of points from   $\widetilde{\Omega}^1_f$ ($\widetilde{\Omega}^1_{f'}$) and homeomorphisms   $\chi_{\sigma}:N_\sigma\to \mathbb{U}$ ($\chi_{\sigma'}:N_{\sigma'}\to \mathbb{U}$) conjugating the restriction of the diffeomorphism  $f^{m_\sigma}$ (${f'}^{m_{\sigma'}}$) on the set   $N_{\sigma}$ ($N_{\sigma'}$) with the diffeomorphism  $b_{\nu}\vert_{_{\mathbb{U}}}$.

Put  $\mathbb{U}_{\tau}=\{(x_{1},...,x_{n})\in
\mathbb{R}^{n}|\ x_{1}^{2}(x_{2}^{2}+...+x_{n}^{2})\leq \tau^2\}$, $\tau\in (0,1]$, recall that we put  $\mathbb{U}_0=\{(x_1,...,x_n)\in \mathbb{R}^n|\ x_1=0\}$ and define  homeomorphisms $\varphi^u_\sigma: \mathbb{U}_0\to \mathbb{U}_0$, $\psi:\mathbb{U}\to \mathbb{U}$  by  $\varphi^u_{\sigma}={\chi}_{\sigma'}\varphi\chi^{-1}_{\sigma}|_{_{\mathbb{U}_0}}$, $\Psi(x_1,x_2,\dots,x_n)=(x_1,\varphi^u_\sigma(x_2,\dots,x_n))$. 

Put $N^\tau_\sigma=\chi^{-1}_\sigma(\mathbb{U}_\tau)$. 
Choose $\tau\in(0,1]$  such that the map $\psi:N^{\tau}_\sigma\to N_{\sigma'}$ defined by  $\psi(x)= \chi_{\sigma'}^{-1}\Psi\chi_\sigma|_{_{N^\tau_\sigma}}$
would be a  well-defined  topological embedding and  $\psi(N^{\tau}_\sigma\setminus W^u_\sigma)\subset \varphi(N^{\tau}_\sigma\setminus W^u_\sigma$).  

Define a topological embedding  $\theta_\sigma:N^{\tau}_\sigma\to N_\sigma$  by  $\theta=\varphi^{-1}\psi$. Since $\theta|_{W^u_\sigma}=id$, it follows from~\cite[Corollary 4.3.2]{GrPo-book} that there exists 
$0<\tau_1<\tau$ and a homeomorphism  $\Theta:N_\sigma\to N_\sigma$ coinciding  with $\theta$ on the set  $N^{\tau_1}_\sigma$ and identical on  $\partial N_\sigma$.

Define homeomorphisms  $h_{\sigma,\sigma'}:N_\sigma\to N_\sigma'$, $h_{O(\sigma),O(\sigma')}:\bigcup \limits_{i=0}^{m_\sigma-1}N_{f^i(\sigma)}\to \bigcup \limits_{i=0}^{m_\sigma-1}N_{{f'}^i(\sigma')}$ by  $h_{\sigma,\sigma'}=\varphi\Theta$, $h_{O(\sigma),O(\sigma')}={f'}^i h_{\sigma,\sigma'}f^{-i}(x)$ for  $x\in N_{f^i(\sigma)}$.

Denote by $H_1:\bigcup \limits_{\sigma\in \Omega_f^1}N_{\sigma} \to \bigcup \limits_{\sigma'\in \Omega_{f'}^1}N_{\sigma'}$ a homeomorphism coinciding for any point  $\sigma\in \Omega_f^1$ with the homeomorphism  $h_{O(\sigma),O(\sigma')}$.

To define the homeomorphism $H_{n-1}:\bigcup \limits_{\sigma\in \Omega_f^{n-1}}N_{\sigma} \to \bigcup \limits_{\sigma'\in \Omega_{f'}^{n-1}}N_{\sigma'}$ use  the similar construction for points from the set $\Omega_f^{n-1}$ using formal replace  $s$ with  $u$ and  $b_{\nu}$ with $b^{-1}_\nu$.

\section{Proof of  Auxiliary Facts} \label{appendix} 
 
\subsection{On Embedding of Families of Closed Curves and Their Tubular Neighborhoods}\label{emb}

In this section we prove Proposition~\ref{nb-equ} after all necessary definitions and auxiliary facts.

A topological space   $X$ is called  {\it $m$-connected} (for $m>0$) if it is non-empty, path-connected and its first $m$  homotopy groups $\pi_i(X)$, $i\in \{1,\dots,m\}$ are trivial. Non-empty  space $X$ can be interpreted as (-1)-connected  and path-connected space  $X$  can be interpreted as 0-connected. 

Let $Q^q, M^n$ be topological manifolds of dimensions $q, n$,  possibly with non-empty boundaries.
 
Topological embeddings $e, e': Q^q\to M^n$ are {\it homotopic} if there is a continuous map $H:Q^q\times [0,1]\to M^n$  ({\it a homotopy}) such that $H|_{Q^q\times \{0\}}=e, H|_{Q^q\times \{1\}}=e'$,    are {\it concordant} if there is a topological embedding $H:Q^q\times [0,1]\to M^n\times [0,1]$ ({\it a concordance}) such that $H(Q^q, 0)=(e(Q^q), 0), H(Q^q, 1)=(e'(Q^q), 1)$ and {\it ambient isotopic with a subset $N\subset M^n$ fixed}, if there is a  family $h_t$ of  homeomorphisms $h_t : M^n\to M^n$, $t\in [0,1]$ ({\it an ambient isotopy}),  such that   $h_0=id$, $h_1{e}=e'$, $h_t|_{_{N}}=id$ for any  $t\in [0,1]$. For $\varepsilon>0$ an ambient isotopy is called {\it an  $\varepsilon$-isotopy} if $d(h_t(x),x)\leq \varepsilon$ for any $x\in M^n, t\in [0,1]$.

 Let $K$ is a simplicial complex in $\mathbb{R}^n$.  A topological space $P$ generated by points of the complex $K$ with the topology induced from  $\mathbb{R}^n$ is called {\it the polyhedron}. 

A complex $K'$ is called {\it a subdivision}  of a complex  $K$ if every simplex of $K'$ belongs to a simplex of $K$. 

 Let $K, L$ be complexes  and $P,Q$ be  polyhedra generated by $K,L$.   A  map  $h:P\to Q$  is called  {\it  piecewise linear} if there exists subdivisions $K', L'$ of $K, L$ correspondingly  such that   $h$ moves each simplex of the complex $K'$ into a simplex of the complex $L'$ (see for example  \cite{RuSa}).

A polyhedron $P$ is called  {\it the  piecewise linear manifold} of dimension  $n$ with boundary if it is a topological manifold with boundary and  for any point  $x\in int\,P$ ($y\in \partial P$)  there is a neighborhood $U_x$ ($U_y$) and a piecewise linear  homeomorphism  $h_x:U_x\to \mathbb{R}^n$ ($h_y:U_y\to \mathbb{R}^n_+=\{(x_1,...,x_n)\subset \mathbb{R}^n|\ x_1\geq 0\}$).

The following important statement  follows from  Theorem~4 of~\cite{Hu}. 

\begin{st}\label{Huds} Suppose that $Q^q, M^n$ are compact  piecewise linear manifolds of dimension  $q, n$ correspondingly,   $Q^q$ is the manifold without boundary,   $M^n$ possibly has a non-empty boundary, $e, e':Q^q\to int~M^n$ are piecewise linear homotopic piecewise linear embeddings, and the following conditions hold:
\begin{enumerate}
\item $q\leq n-3$;
\item $Q^q$ is $(2q-n+1)$-connected;
\item $M^n$ is $(2q-n+2)$-connected.
\end{enumerate}

Then $e,e'$ are piecewise linear ambient isotopic with $\partial M^n$ fixed.   
\end{st}

To reduce the problem of embedding of families of curves to the statement~\ref{Huds} we use the following results.

The combination of Theorems 1 and 4 of R. Miller's paper~\cite{Mi}, see also A. Chernavskii paper~\cite{Ch}, gives the next statement.

\begin{st}\label{Mi-Ch} Let $Q^q, M^n$ be compact  piecewise linear manifolds of dimension  $q, n$ correspondingly, $q\leq n-3$, $\varepsilon>0$, and $e: N^k\to M^n$ be a locally flat embedding. Then there is an ambient $\varepsilon$-isotopy of $M^n$ connecting $e$ with piecewise embedding.   
\end{st}

G. Weller's result~\cite[Theorem 3]{We} gives the following. 

\begin{st}\label{Wel} Let $Q^q, M^n$ be topological manifolds of dimension  $q, n$ correspondingly, $Q^q$ is compact,   $e, e':Q^q\to int~M^n$ are  homotopic embeddings, and the following conditions hold:
\begin{enumerate}
\item $q\leq n-3$;
\item $Q^q$ is $(2q-n+1)$-connected;
\item $M^n$ is $(2q-n+2)$-connected.
\end{enumerate}

Then $e,e'$ are locally flat concordant, that is there is a locally flat   embedding $H:Q^q\times [0,1]\to M^n\times [0,1]$  such that $H(Q^q, 0)=(e(Q^q), 0), H(Q^q, 1)=(e'(Q^q), 1)$.
\end{st}


\begin{Prop}\label{triv-000} Suppose that $\mathbb{P}^{n-1}$ is either $\mathbb{S}^{n-1}$ or  $\mathbb{B}^{n-1}$,  $\{\beta_i\}, \{\beta'_i\}\subset int\,\mathbb{P}^{n-1}\times \mathbb{S}^1$ are two families of pairwise disjoint knots such that knots $\beta_i,\beta'_i$ are homotopic, $i\in \{1,\dots,k\}$, $n\geq 4$.  Then there exists a homeomorphism   $H: \mathbb{P}^{n-1}\times \mathbb{S}^1 \to \mathbb{P}^{n-1}\times \mathbb{S}^1$ such that    $H(\beta_i)=\beta'_i$ for any $i\in \{1,\dots,k\}$. Moreover, if $\mathbb{P}^{n-1}=\mathbb{B}^{n-1}$ then    $H|_{_{\partial\,\mathbb{P}^{n-1}\times \mathbb{S}^1}}=id$. 
\end{Prop}
\begin{demo} Denote by $\Delta^n$ the simplex of dimension $n$ and by $\partial\,\Delta^n$ its boundary. Then $\mathbb{P}^{n-1}\times \mathbb{S}^1$ is homeomorphic either to $\partial\, \Delta^{n}\times \partial \Delta^2$ or to $\Delta^{n-1}\times \partial \Delta^2$ and without loss of generality we will identify $\mathbb{P}^{n-1}\times \mathbb{S}^1$  with one of  these piecewise linear objects.

  It follows from Statement~\ref{Mi-Ch} that there exist  homeomorphisms  $g,g':\mathbb{P}^{n-1}\times \mathbb{S}^1\to \mathbb{P}^{n-1}\times \mathbb{S}^1$ such that for any      $i\in \{1,...,k\}$  the sets   $g(\beta_i), g'(\beta'_i)$   are  subpolyhedra.  It is sufficient to construct a homeomorphism  $H: \mathbb{P}^{n-1}\times \mathbb{S}^1 \to \mathbb{P}^{n-1}\times \mathbb{S}^1$ such that     $H(g(\beta_i))= g'(\beta'_i)$ for any $i\in \{1,\dots,k\}$ and    $H|_{_{\partial\,\mathbb{P}^{n-1}\times \mathbb{S}^1}}=id$. Then the homeomorphism ${g'}^{-1}Hg$ will be the desired map.  
So without loss of generality assume that knots  $\beta_i, \beta'_i$  are sub-polyhedra for  $i\in \{1,\dots,k\}$.

By assumption     piecewise linear embeddings  ${e}_i:\partial \Delta^2\to  \mathbb{P}^{n-1}\times \mathbb{S}^1$, $e'_i: \partial \Delta^2\to \mathbb{P}^{n-1}\times \mathbb{S}^1$  such that   ${e}_i(\partial \Delta^2)=\beta_1$, ${e}'_i(\partial \Delta^2)=\beta'_i$ are homotopic. By Statement~\ref{Wel}  there is a locally flat   embedding $\Psi:\partial\Delta^2\times [0,1]\to \mathbb{P}^{n-1}\times \mathbb{S}^1\times [0,1]$  such that $\Psi(\partial\Delta^2, 0)=(e(\partial\Delta^2), 0), \Psi(\partial\Delta^2, 1)=(e'(\partial\Delta^2), 1)$. By Statement~\ref{Mi-Ch} there is a homeomorphism $\Phi: \mathbb{P}^{n-1}\times \mathbb{S}^1\times [0,1]\to \mathbb{P}^{n-1}\times \mathbb{S}^1\times [0,1]$ such that superposition  $\Phi\Psi$ is piecewise linear embedding. Denote by $pr:\mathbb{P}^{n-1}\times \mathbb{S}^1\times [0,1]\to \mathbb{P}^{n-1}\times \mathbb{S}^1\times \{0\}$ the projection given by $pr(x,t)=(x,0)$, $x\in \mathbb{P}^{n-1}\times \mathbb{S}^1$, $t\in [0,1]$. Then superposition $pr\Phi\Psi$ is the piecewise linear homotopy.   

Now  construct the homeomorphism $H$  by induction on $i$. 
By   Statement~\ref{Huds}, there exists a   piecewise linear homeomorphism   $H_1: \mathbb{P}^{n-1}\times \mathbb{S}^1\to \mathbb{P}^{n-1}\times \mathbb{S}^1$ that   maps $\beta_1$ to $\beta'_1$. 

Suppose that  for some  $j<k$  there exists a homeomorphism $H_j:\mathbb{P}^{n-1}\times \mathbb{S}^1\to \mathbb{P}^{n-1}\times \mathbb{S}^1$ such that $H_j(\beta_i)=\beta'_i$ for $i\in \{1,\dots,j\}$,  $H_j|_{\partial \mathbb{P}^{n-1}\times \mathbb{S}^1}=id$. Construct a homeomorphism  $H_{j+1}:\mathbb{P}^{n-1}\times \mathbb{S}^1\to \mathbb{P}^{n-1}\times \mathbb{S}^1$ such that $H_{j+1}(\beta_i)=\beta'_i$ for $i\in \{1,\dots,j+1\}$,  $H_{j+1}|_{\partial \mathbb{P}^{n-1}\times \mathbb{S}^1}=id$.    

  Preserve notation  $\beta_i$ for images of knots $\beta_i$, $i\in \{1,\dots,k\}$, by means of homeomorphism $H_j$ (then the knots $\beta_i, \beta'_i$  coincide for $i\in \{1,\dots,j\}$) and denote by $e_{j+1}, e'_{j+1}:\partial \Delta^2\to \mathbb{P}^{n-1}\times \mathbb{S}^1$ piecewise linear embeddings such that $e_{j+1}(\partial \Delta^2)=\beta_{j+1}, e'_{j+1}(\partial \Delta^2)=\beta'_{j+1}$. 
	
	Show that $e_{j+1}, e'_{j+1}$ are pieciwise  linear homotopic in $\mathbb{P}^{n-1}\times \mathbb{S}^1\setminus \bigcup\limits_{i=1}^{j}\beta_i$.  	It follows from   Statement~\ref{Huds} that  there exists a   family of piecewise linear homeomorphism   $h^{j+1}_t : \mathbb{P}^{n-1}\times \mathbb{S}^1\to \mathbb{P}^{n-1}\times \mathbb{S}^1$, $t\in [0,1]$,  such that   $h^{j+1}_0=id$, $h^{j+1}_1{e_{j+1}}=e'_{j+1}$, $h^{j+1}_t|_{_{\partial \mathbb{P}^{n-1}\times \mathbb{S}^1}}=id$ for any  $t\in [0,1]$. This family defines a piecewise linear homeomorphism  $\overline{H}^{j+1}:\mathbb{S}^1\times [0,1]\to \mathbb{P}^{n-1}\times \mathbb{S}^1\times [0,1]$ by $\overline{H}^{j+1}(x,t)=(h^{j+1}_t(x),t)$. It follows from  general position Theorem (see, for example, Theorem 5.3 of Chapter 5 of \cite{RuSa}) that  there   exists  a piecewise linear homeomorphism $\Phi: \mathbb{P}^{n-1}\times \mathbb{S}^1\times [0,1]\to \mathbb{P}^{n-1}\times \mathbb{S}^1\times [0,1]$ leaving the set $\partial\mathbb{P}^{n-1}\times \mathbb{S}^1\times [0,1]\cup \mathbb{P}^{n-1}\times \mathbb{S}^1\times \{0,1\}$ fixed and  mapping the cylinder $\overline{H}^{j+1}(\mathbb{S}^1\times [0,1])$ to one that have no intersection with all cylinders $\beta_{i}\times [0,1]$, $i\in \{1,\dots,j\}$.  Then it possible to choose regular tubular neighborhoods  $N_{\beta'_1},\dots, N_{\beta'_j}$ of knots $\beta'_1,\dots,\beta'_j$ in $\mathbb{P}^{n-1}\times \mathbb{S}^1$ such that $\Phi(\overline{H}^{j+1}(\mathbb{S}^1\times [0,1]))\cap \bigcup \limits_{i=1}^{j} N_{\beta'_j}\times [0,1]=\emptyset$. So, the map $\Phi\overline{H}^{j+1}:\mathbb{S}^1\times [0,1]\to (\mathbb{P}^{n-1}\times \mathbb{S}^1\setminus \bigcup \limits_{i=1}^{j} int\,N_{\beta'_j})\times [0,1]$ is a piecewise concordance and composition of it with projection $pr:\mathbb{P}^{n-1}\times\mathbb{S}^1\times [0,1]\to \mathbb{P}^{n-1}\times\mathbb{S}^1$ give the desire homotopy.  Applying Statement~\ref{Huds} once more one can obtain the piecewise linear homeomorphisms   $\tilde{H}^{j+1}: \mathbb{P}^{n-1}\times \mathbb{S}^1\setminus \bigcup \limits_{i=1}^{j}int\,N_{\beta'_j} \to \mathbb{P}^{n-1}\times \mathbb{S}^1\setminus \bigcup \limits_{i=1}^{j}int\,N_{\beta'_j}$, $t\in [0,1]$,  such that    $\tilde{H}^{j+1}{e_{j+1}}=e'_{j+1}$, $\tilde{H}^{j+1}|_{_{\partial (\mathbb{P}^{n-1}\times \mathbb{S}^1\setminus \bigcup \limits_{i=1}^{j}int\,N_{\beta'_j})}}=id$ for any  $t\in [0,1]$. Now  the desire homeomorphism $H^{j+1}: \mathbb{P}^{n-1}\times \mathbb{S}^1\to \mathbb{P}^{n-1}\times \mathbb{S}^1$ is defined by 
$H^{j+1}(x)=\tilde{H}^{j+1}(x)$ for  $x\in \mathbb{P}^{n-1}\times \mathbb{S}^1\setminus \bigcup \limits_{i=1}^{j}int\,N_{\beta'_j}$ and  $H^{j+1}(x)=x$ for $x\in \bigcup \limits_{i=1}^{j}int\,N_{\beta'_j}$. 

 \end{demo}

The following Statement~\ref{prod}  is proved in~\cite[Lemma~2.1]{GrGuMe10}. 

\begin{st}\label{prod} Let $h:\mathbb{B}^{n-1}\times \mathbb{S}^1\to  int~\mathbb{B}^{n-1}\times \mathbb{S}^1$ be a topological embedding  such that  $h(\{O\}\times \mathbb{S}^1)=\{O\}\times \mathbb{S}^1$.  Then a manifold 
 $\mathbb{B}^{n-1}\times \mathbb{S}^1\setminus int~h(\mathbb{B}^{n-1}\times \mathbb{S}^1)$ is homeomorphic to the direct product  $\mathbb{S}^{n-2}\times \mathbb{S}^1\times [0,1]$.
\end{st}

The proof of the following statement  one can find in~\cite[Lemma 2]{Le}.
\begin{st}\label{plus-prod} Suppose that  $Y$ is a topological manifold with boundary, $X$ is a connected component of its boundary, $Y_1$ is a manifold homeomorphic to $X\times [0,1]$, and   $Y\cap Y_1=X$.  Then a manifold $Y\cup Y_1$ is homeomorphic to $Y$. Moreover,   if the manifold $Y$ is homeomorphic to the direct product  $X\times [0,1]$ then there exists a homeomorphism  $h:   X\times [0,1]\to Y\cup Y_1$  such that   $h(X\times \{\frac{1}{2}\})=X$.
 \end{st}

\medskip
{\bf Proof of Proposition~\ref{nb-equ}.} 

Suppose that $\mathbb{P}^{n-1}$ is either the ball  $\mathbb{B}^{n-1}$ or the sphere  $\mathbb{S}^{n-1}$,   $\{\beta_i\}, \{\beta'_i\}\subset int\,\mathbb{P}^{n-1}\times \mathbb{S}^1$   are two families of knots such that knots $\beta_i,\beta'_i$ are homotopic and $\{N_{\beta_i}\}, \{N_{\beta'_i}\}\subset \mathbb{P}^{n-1}\times \mathbb{S}^1$ are   pairwise disjoint neighborhoods of knots $\{\beta_i\}, \{\beta'_i\}$.  Let us prove that  there exist a homeomorphism  $h:\mathbb{P}^{n-1}\times \mathbb{S}^1\to \mathbb{P}^{n-1}\times \mathbb{S}^1$  such that $h(\beta_i)=\beta'_i, h({N}_{\beta_i})={N}_{\beta'_i}, i\in \{1,...,k\}$ and  $h|_{_{\partial\,\mathbb{P}^{n-1}\times \mathbb{S}^1}}=id$.

  By  Proposition~\ref{triv-000}, there exists a homeomorphism   $h_0: \mathbb{P}^{n-1}\times \mathbb{S}^1\to \mathbb{P}^{n-1}\times \mathbb{S}^1$ such that  $h_0(\beta_i)=\beta'_i$, $h_0|_{_{\partial\, \mathbb{P}^{n-1}\times \mathbb{S}^1}}=id$. Put  $\tilde{N}_i=h_0(N_{\beta_i})$. For $r\in (0,1)$ put $\mathbb{B}^{n-1}_r=\{(x_1,\dots,x_{n-1})\in \mathbb{R}^{n-1}|\ x_1^2+\dots+x_{n-1}^2\leq r^2\}$.  It follows from~\cite{Br62} that  there exist topological embeddings  $\tilde{e}_i:\mathbb{B}^{n-1}\times\mathbb{S}^1 \to  int\,\mathbb{P}^{n-1}\times \mathbb{S}^1$ such that  $\tilde{e}_i(\{O\}\times\mathbb{S}^1)=\beta_i'$,   $\tilde{e}_i(\mathbb{B}_r^{n-1}\times\mathbb{S}^1)=\tilde{N}_i$ for some $r\in (0,1)$, $\tilde{e}_i(\mathbb{B}^{n-1}\times\mathbb{S}^1) \cap \tilde{e}_j(\mathbb{B}^{n-1}\times\mathbb{S}^1)=\emptyset$ for $i\neq j, i,j\in \{1,...,k\}$. Put $U_i=\tilde{e}_i(\mathbb{B}^{n-1}\times\mathbb{S}^1)$. 

Denote by $e'_i:\mathbb{B}^{n-1}\times \mathbb{S}^1\to \mathbb{P}^{n-1}\times \mathbb{S}^1$ a topological embedding such that $e'_i(\mathbb{B}^{n-1}\times\mathbb{S}^1)=N_{\beta_i'}$ and choose $r_0, r_1$ such that $0<r_0<r_1<1$ and $e'_i(\mathbb{B}^{n-1}_{r_1}\times \mathbb{S}^1)\subset \tilde{N}_i$. Put $N'_{0,i}=e'_i(\mathbb{B}^{n-1}_{r_0}\times \mathbb{S}^1), N'_{1,i}=e'_i(\mathbb{B}^{n-1}_{r_1}\times \mathbb{S}^1)$.

By  Statement~\ref{prod} the  set  $\tilde{N}_i\setminus N'_{1,i}$ is homeomorphic to the direct product $\mathbb{S}^{n-2}\times \mathbb{S}^1\times [0,1]$, $i\in \{1,\dots,k\}$. By  Statement~\ref{plus-prod}, there exists a homeomorphism  $g_i:\mathbb{S}^{n-2}\times \mathbb{S}^1\times [0,1]\to U_i\setminus int\,N'_{0,i}$ and $t_1,t_2\subset (0,1)$ such that $g_i(\mathbb{S}^{n-2}\times \mathbb{S}^1\times \{t_1\})=\partial\, \tilde{N}_i$, $g_i(\mathbb{S}^{n-2}\times \mathbb{S}^1\times \{t_2\})=\partial\, N'_{1,i}$.   Let $\xi:[0,1]\to [0,1]$ be a homeomorphism that is identity on the ends of the interval  $[0,1]$ and such that $\xi(t_1)=t_2$.  Define a homeomorphism  $\tilde{g}_i:\mathbb{S}^{n-2}\times \mathbb{S}^1\times [0,1]\to  \mathbb{S}^{n-2}\times \mathbb{S}^1\times [0,1]$  by  $\tilde{g}_i(x,t)=(x,\xi(t))$ and  a homeomorphism $h_i:\mathbb{P}^{n-1}\times \mathbb{S}^1\to \mathbb{P}^{n-1}\times \mathbb{S}^1$  by

$$h_i(x)=\begin{cases}
g_i(\tilde{g}_i(g^{-1}_i(x))),~x\in U_i\setminus int\,N'_{0,i};
\cr x,~x\in (\mathbb{P}^{n-1}\times \mathbb{S}^1\setminus U_i).\end{cases}$$

The superposition  $\eta=h_k\cdots h_1h_0$ maps every knot  $\beta_i$ into  the knot  $\beta'_i$,   the  neighborhood $N_{\beta_i}$ into the set    $N'_{1,i}\subset N_{\beta'_i}$, and keeps  the set  $\partial\, \mathbb{P}^{n-1}\times \mathbb{S}^1$ fixed.  Since the set $N_{\beta'_i}\setminus int\,N'_{1,i}$ is homeomorphic to the direct product $\mathbb{S}^{n-2}\times \mathbb{S}^1\times [0,1]$, it is possible to  apply the  described construction  once more and  get the desire homeomorphism. The proof is complete.

Let $M^n$ be   homeomorphic to   $\mathbb{S}^{n-1}\times \mathbb{S}^1$. The knot   $\beta\in int\,M^n$ is called {\it essential} if a  homomorphism  $i_*:\pi_1(\beta)\to \pi_1(M^n)$ induces by an inclusion  $i:\beta \to {M}^n$  is an isomorphism.  

The next statement follows immediately from Proposition~\ref{nb-equ}.
 
\begin{Cor}\label{blue}
 If $\beta\subset M^n$ is an essential knot  and $N_\beta$ is its  tubular neighborhood  then the manifold  $M^n\setminus int~N_\beta$ is homeomorphic to the direct product  $\mathbb{B}^{n-1}\times \mathbb{S}^1$.
\end{Cor}

One more corollary was stated in Subsection~\ref{start1}, let us prove it.

We will use the  following statement proved in~\cite[Theorem 2]{Max}.
 
\begin{st}\label{maxx} Let  $\psi:\mathbb{S}^{n-2}\times \mathbb{S}^1\to \mathbb{S}^{n-2}\times \mathbb{S}^1$  be an arbitrary homeomorphism.  Then there exists a homeomorphism    $\Psi:\mathbb{B}^{n-1}\times \mathbb{S}^1\to \mathbb{B}^{n-1}\times \mathbb{S}^1$  such that   $\Psi|_{_{\mathbb{S}^{n-2}\times \mathbb{S}^1}}~=~\psi|_{_{\mathbb{S}^{n-2}\times \mathbb{S}^1}}$.
\end{st}

\medskip
{\bf Proof of Corollary~\ref{contin}.}

Let $\tilde{h}:\partial{\mathbb{B}^{n-1}\times \mathbb{S}^1}\to \partial{\mathbb{B}^{n-1}\times \mathbb{S}^1}$ be a homeomorphism. Prove that there  exists  a homeomorphism  $h:\mathbb{B}^{n-1}\times \mathbb{S}^1\to \mathbb{B}^{n-1}\times \mathbb{S}^1$  such that $h(\beta_i)=\{x_i\}\times \mathbb{S}^1, h({N}_{\beta_i})={N}_{\beta'_i}, i\in \{1,...,k\}$, and  $h|_{_{\partial\,\mathbb{B}^{n-1}\times \mathbb{S}^1}}=\tilde{h}$.  Due to Statement~\ref{maxx} the homeomorphism $\tilde{h}$ can be extended  to a homeomorphism $h_0: \mathbb{B}^{n-1}\times \mathbb{S}^1\to \mathbb{B}^{n-1}\times \mathbb{S}^1$. Due to Proposition~\ref{nb-equ} there exists a homeomorphism $h_1:\mathbb{B}^{n-1}\times \mathbb{S}^1\to \mathbb{B}^{n-1}\times \mathbb{S}^1$  such that $h_1(h_0(\beta_i))=\beta'_i, h_1(h_0({N}_{\beta_i}))={N}_{\beta'_i}, i\in \{1,...,k\}$ and  $h_1|_{_{\partial\,\mathbb{B}^{n-1}\times \mathbb{S}^1}}=id$. Then the superposition $h_1h_0$ is the desire homeomorphism $h$.

\subsection{Characteristic Space and Embedding of Separatrices of Dimension $(n-1)$}\label{orbits1}

Here we prove Lemmas~\ref{svo1},\ref{svo}.  If all saddle points of the diffeomorphism $f$ are fixed and has positive orientation type then Lemmas  follows from~\cite[Lemm 3.1]{GrGuPo-emb}.  Prove them for general case.  The main tool of  the proof is a surgery along   knots that,  in contrast with case $n=3$, does not change the topology of a  manifold.  

Let  $M^n$ be a topological manifold of dimension $n\geq 4$ (possibly with non-empty boundary),  $\beta\in int\,M^n$ be a knot  and $N_\beta\subset int\, M^n$ be its tubular neighborhood. Glue manifolds $M^n\setminus int\, N_{\beta}$ and  $\mathbb{B}^{n-1}\times \mathbb{S}^1$ by means of an arbitrary reversing the  natural orientation homeomorphism $\varphi:\partial N_\beta\to \mathbb{S}^{n-2}\times \mathbb{S}^1$ and denote the obtained manifolds by $Q^n$. We say that $Q^n$ is obtained from $M^n$ by {\it the surgery along the knot $\beta$}.

\begin{Prop}\label{ext}  $Q^n$ is homeomorphic to  $M^n$. 
\end{Prop}
\begin{demo} Put  $N'=M^n\setminus int~N_\beta$ then  $Q^n=N'\cup_{\varphi}\mathbb{B}^{n-1}\times \mathbb{S}^1$ and for any  subset  $X\subset  N'\cup \mathbb{B}^{n-1}\times \mathbb{S}^1$ a projection   $\pi: X\to Q^n$ is defined.  

 Put $\psi= \varphi^{-1}\pi^{-1}|_{\pi(\mathbb{S}^{n-2}\times \mathbb{S}^1)}$. Due to Statement~\ref{maxx} the homeomorphism   $\psi$ can be extended up to a homeomorphism  $\Psi: \pi(\mathbb{B}^{n-1}\times \mathbb{S}^1)\to N_\beta$.  Then a map  $H:Q^n\to M^n$, defined by   $H(x)=\pi^{-1}(x)=x$ for  $x\in \pi(int\,N')$ and by  $H(x)=\Psi(x)$ for  $x\in \pi(\mathbb{B}^{n-1}\times \mathbb{S}^1)$ is the  desire homeomorphism. 
\end{demo}

Recall that we represent the sphere   $S^n$  as the union of pairwise disjoint sets  $A_f=(\bigcup\limits_{\sigma\in \Omega^{1}_f}{W^u_\sigma})\cup \Omega^{0}_f,\, R_f=(\bigcup\limits_{\sigma\in \Omega^{n-1}_f}{W^s_\sigma})\cup  \Omega^{n}_f,\,V_f=S^n\setminus(A_f\cup R_f)$, denoted by  $\widehat V_f=V_f/f$  the orbit space of the action of  $f$ on $V_f$,  by   $p_{_f}:V_f\to \widehat V_f$ the natural projection and   introduced the homomorphism $\eta_{_{f}}:\pi_{1}(\widehat V_f)\to \mathbb{Z}$. 

Suppose that the set $\Omega^1_f$ is non-empty (in opposite case consider $f^{-1}$). Then $A_f$ can be represented as a graph embedded in $S^n$ whose vertices are sink periodic points and edges are one-dimensional unstable manifolds of saddle periodic points. Since the closure of all stable  separatrix of dimension $(n-1)$  cuts the ambient sphere  $S^n$ into  disjointed union of  domain each of which contains exactly one sink periodic point then  $|\Omega^0_f|=|\Omega^1_f|+1$.  It follows from~\cite{GrPoZh} that  the sets  $A_f$ is connected.  So, due to Proposition~\ref{tr}, $A_f$ is a tree  and it is possible to define {\it ranks of sink periodic points} as   ranks of the vertices  of $A_f$.  


Construct a series of   attractors $A_0,A_1,\dots,A_r$, where  $A_0=\bigcup\limits_{\omega\in \Omega^0_f}\omega$, $A_r=A_f$, $A_{i}$ consists of all  vertex of $A_f$ of ranks less of equal to $i$ and joining them edges of $A_f$.  Denote by $V_i$ the union of stable manifolds of all periodic points belonging to $A_i$,  put $\widehat{V}_i=V_i/_f$ and denote by $p_i:V_i\to \widehat{V}_i$ the natural projection. It follows from the definition that the number of connected components equals to the number of $f-$invariant components of $A_i$.

Lemma~\ref{svo1} immediately follows from the next proposition.

\begin{Prop}\label{ind}
$\widehat{V}_i$ is a union of manifolds homeomorphic to $\mathbb{S}^{n-1}\times \mathbb{S}^1$.
\end{Prop}
\begin{demo}  Prove the proposition by induction on $i$.   
		For a point $\omega\in \Omega^0_f$ of period $m_\omega$ put   $V^s_{\mathcal{O}_\omega}=\bigcup \limits_{i=0}^{m_\omega-1} f( W^s_\omega\setminus \omega)$ and  $\widehat{V}^s_{\mathcal{O}_\omega}=V^s_\omega/f$. It  follows  from hyperbolicity of the point $\omega$ and~\cite[Theorem~5.5]{Ko}  (see also~\cite[Propositions 1.2.3 и 1.2.4]{BoGrPo05}) that $\widehat{V}^s_{\mathcal{O}_\omega}$ is homeomorphic to  $\mathbb{S}^{n-1}\times \mathbb{S}^1$. So $\widehat{V}_0$ is the union of manifolds homeomorphic to $\mathbb{S}^{n-1}\times \mathbb{S}^1$.  

Suppose the statement is proven for $k=i$ and prove it for $i+1$.

For a point  $\sigma\in \Omega^{1}_f$ of period $m_\sigma$   denote by $\omega_{-}, \omega_{+}$ the sink points belonging to $cl\,W^u_\sigma$ and by $l^u_{\sigma,-}, l^u_{\sigma,+}$ unstable separatrices of $\sigma$ such that $l^u_{\sigma,-}\subset W^u_{\omega_{+}}, l^u_{\sigma,+}\subset W^u_{\omega_{+}}$.  Put $\mathcal{O}_{\sigma}=\bigcup\limits_{i=0}^{m_\sigma-1}f^{i}(\sigma)$, $l^u_{\mathcal{O}_\sigma,-}= \bigcup\limits_{i=0}^{m_\sigma-1}f^{i}(l^u_{\sigma,-})$, $l^u_{\mathcal{O}_\sigma,+}= \bigcup\limits_{i=0}^{m_\sigma-1}f^{i}(l^u_{\sigma,+})$.   Let $N_{\sigma}$ be the linearizing neighborhood of the point $\sigma$ and $N^s_{\sigma}=N_\sigma\setminus W^s_\sigma$, $N^u_{\sigma}=N_\sigma\setminus W^u_\sigma$.  Denote by $N^u_{\sigma,-}, N^u_{\sigma,+}$ connected components of the set $N^u_{\sigma}$ containing separatrices $l_{\sigma,-},l_{\sigma,+}$ correspondingly, put $\widehat{N}^u_{\sigma}=N^u_{\sigma}/f$, $\widehat{N}^s_{\sigma}=N^s_{\sigma}/f$,   denote by $p_{\widehat{N}^u_{\sigma}}: N^u_{\sigma}\to \widehat{N}^u_{\sigma}, p_{\widehat{N}^s_{\sigma}}: N^s_{\sigma}\to \widehat{N}^s_{\sigma}$  the natural projections and define a homeomorphism $\varphi: \partial\,\widehat{N}^u_{\sigma}\to \widehat{N}^s_{\sigma}$ by $\varphi=p_{\widehat{N}^s_{\sigma}}p^{-1}_{\widehat{N}^u_{\sigma}}$.

Suppose that $W^u_\sigma\subset A_{i+1}\setminus A_i$ and $\omega_-\subset A_i$. Two cases are possible: 1)  $rank(\omega_{-})<rank (\omega_{+})$, then  $\omega_{+}\subset A_{i+1}\setminus A_i$; 2) $rank(\omega_{-})=rank(\omega_{+})$, then   $\omega_{+}\subset A_i$.

Consider case   1). It follows from Proposition~\ref{dima-c1}  that  $m_\sigma=m_{\omega_{-}}$, moreover, the period of the connected component $V^{\sigma}_i$ of $V_i$ having non-empty intersection with the set  $W^s_{\omega_{-}}$ also equals $m_\sigma$.  Then the set $\hat{l}^u_{\mathcal{O}_{\sigma},-}=p_{i}(l^u_{\mathcal{O}_\sigma,-})$ is an essential knot in   $\widehat{V}^{\sigma}_i=V^{\sigma}_i/_f$ and  due to  Corollary~\ref{blue}  the set $\widehat{V}^{\sigma}_{i}\setminus int\, \widehat{N}^u_{\sigma,-}$ is homeomorphic to $\mathbb{B}^{n-1}\times \mathbb{S}^1$.

Denote by $p_{\widehat{V}_{\mathcal{O}_{\omega_{+}}}}: {V}_{\mathcal{O}_{\omega_{+}}}\to  \widehat{V}_{\mathcal{O}_{\omega_{+}}}$ the natural projection. Then the  set $\hat{l}^u_{\mathcal{O}_{{\sigma},+}}=p_{\widehat{V}_{\mathcal{O}_{\omega_{+}}}}(l^u_{\mathcal{O}_{\sigma,+}})$ is a knot in $\widehat{V}_{\omega_{\sigma,+}}$ and  $\widehat{N}^u_{\sigma,+}=p_{\widehat{V}_{\mathcal{O}_{\omega_{+}}}}(N^u_{\sigma,+})$ is its tubular neighborhood. 

Put $V^\sigma_{i+1}=V^\sigma_i\cup {V}_{\mathcal{O}_{\omega_+}}\cup W^s_{\mathcal{O}_\sigma} \setminus W^u_{\mathcal{O}_\sigma}=V^\sigma_i\cup {V}_{\mathcal{O}_{\omega_+}}\cup N^s_{\mathcal{O}_\sigma} \setminus int\, N^u_{\mathcal{O}_\sigma}=(V^\sigma_i\setminus int\, N^u_{\mathcal{O}_\sigma,-})\cup N^s_{\mathcal{O}_\sigma}  \cup ({V}_{\mathcal{O}_{\omega_+}} \setminus int\, N^u_{\mathcal{O}_{\sigma,+}})$. 
Then $\widehat{V}^\sigma_{i+1}=V^\sigma_{i+1}/f=((\widehat{V}_{\mathcal{O}_{\omega_-}}\setminus int\,\widehat{N}^u_{\sigma,-})\cup_{\varphi_-}\widehat{N}^s_\sigma)\cup_{\varphi_+} (\widehat{V}^\sigma\setminus int\,\widehat{N}^u_{\sigma,-})$, where $\varphi_-=\varphi|_{\partial\widehat{N}^u_{\sigma,-}}$, $\varphi_+=(\varphi|_{\partial\widehat{N}^u_{\sigma,+}})^{-1}$.
 
Due to Proposition~\ref{st-nbh} the  manifold  $\widehat{N}^s_\sigma$ is homeomorphic to $\mathbb{S}^{n-2}\times \mathbb{S}^1\times [-1,1]$. 
Due to Statement~\ref{plus-prod} $(\widehat{V}_{\mathcal{O}_{\omega_+}}\setminus int\,\widehat{N}^u_{\sigma,-})\cup_{\varphi_-}\widehat{N}^s_\sigma$ is homeomorphic to $\widehat{V}_{\mathcal{O}_{\omega_-}}\setminus int\,\widehat{N}^u_{\sigma,-}$, so, is homeomorphic to $\mathbb{B}^{n-1}\times \mathbb{S}^1$.   Then $\widehat{V}^\sigma_{i+1}$ is obtained from $\widehat{V}_{\mathcal{O}_{\omega_+}}\setminus int\,\widehat{N}^u_{\sigma,+}$ 
by surgery along knot $\hat{l}^u_{\sigma,+}$ and, due to Proposition~\ref{ext}, 
is homeomorphic to $\mathbb{S}^{n-1}\times \mathbb{S}^1$. To get $\widehat{V}_{i+1}$ join $\widehat{V}^\sigma_{i+1}$ to the union  $\bigcup\limits_{rank{\omega}=i+1}\widehat{V}_{\mathcal{O}_\omega}$ and then repeat similar procedure for all saddles $\tilde{\sigma}$ such that $W^u_{\tilde{\sigma}}\subset A_{i+1}\setminus A_i$ and  $A_i\cap cl\,W^u_{\tilde{\sigma}}\neq \emptyset$.  At every step one gets the union of manifolds homeomorphic to $\mathbb{S}^{n-1}\times \mathbb{S}^1$   and after finite number of steps  one gets either to one sink of maximal rank (so the last manifold is connected)  or to  case 2).

In case 2) there are two possibility:   the point $\sigma$ has  either positive or  negative   orientation type. If $\sigma$ has positive orientation type then $m_\sigma=m_{\omega,+}=m_{\omega,-}=1$ and the  manifold $\widehat{V}_i=\widehat{V}_{r-1}$  is a union of two connected components homeomorphic to $\mathbb{S}^{n-1}\times \mathbb{S}^1$ each of which contains one of the knots $\hat{l}^u_{\sigma,+}, \hat{l}^u_{\sigma,-}$. Then  surgery and  arguments similar to ones above prove  that $\widehat{V}_{r}$ is homeomorphic to $\mathbb{S}^{n-1}\times \mathbb{S}^1$.

 If $\sigma$ has negative  orientation type then $m_\sigma=1$,  points $\omega,+, \omega,-$ generate a 2-periodic orbit, and $V_i$ consists of two 2-periodic connected components $D$, $f(D)$. Consider a map $g=f^2$, put $\widehat{V}_g=V_f/_g$ and denote by $p_g:V_f\to \widehat{V}_g$ the natural projection. It follows from  arguments above  that factor-space $V_f/_{f^2}$ is homeomorphic to $\mathbb{S}^{n-1}\times \mathbb{S}^1$. Put $\tau=p_gfp_g^{-1}$. Then $\tau^2=id$, so $\tau$ is involution and $\widehat{V}_f=\widehat{V}_g/_\tau$. It follows from~\cite{Kw2011} that $\widehat{V}_f$   is homeomorphic to one of the following manifolds:  the direct product $\mathbb{S}^{n-1}\times \mathbb {S}^1$, non-oriented fiber bundle  $\mathbb{S}^{n-1} \tilde{\times}\mathbb{S}^1$  over the circle with the fiber  $\mathbb{S}^{n-1}$, the direct product  $\mathbb{S}^1\times \mathbb{R}P^{n-1}$ or  to the connected sum   $\mathbb{R}P^{n}\times \mathbb{R}P^{n}$.

 $\widehat{V}_f$ is covered by $\widehat{V}_g$ and, consequently, by $\mathbb{S}^{n-1}\times \mathbb{R}$, which is the universal cover for  $\widehat{V}_f$. Then due to~\cite[Corollary 19.4]{Ko} the fundamental group  $\pi_1(\widehat V_{f})$ is isomorphic to the group  $\{f^n\}$ and, consequently, to the group $\mathbb{Z}$.  So, $\widehat{V}_f$ cannot be homeomorphic to the direct product  $\mathbb{S}^1\times \mathbb{R}P^{n-1}$ or  to  $\mathbb{R}P^{n}\times \mathbb{R}P^{n}$.  Since   $f$ is orientation preserving then the orbit space $\widehat{V}_f$ is orientable, so it is homeomorphic to $\mathbb{S}^{n-1}\times \mathbb{S}^1$. 
\end{demo}

\medskip
{\bf Proof of Lemma~\ref{svo}.}
\medskip

Prove that $\Delta_{f,i}=\hat{\xi}(\widehat{\Gamma}_f\setminus \widehat{\Gamma}_{f,i+1})$ is  the union of $\mathbb{S}^{n-1}\times \mathbb{S}^1$ for an arbitrary  $i\in \{0,\dots,s\}$.   

If graph $\widehat{\Gamma}_f$ has no loop then for any $i\in \{0,\dots,s\}$ a   connected component   of the set $\widehat{\Gamma}_f\setminus \widehat{\Gamma}_ {f,i+1}$ consist of an edge $e_\sigma$, one of incident to $e_\sigma$ vertex $v$ and all paths of the graph $\widehat{\Gamma}_f$ connecting the vertex $v_\sigma$ with leaves and crossing the vertices in order of decreasing of ranks. Denote this component by $T^\sigma$ and put  $\Delta^\sigma_{f,i}=\widehat{\xi}(T^\sigma)$. Prove that  $\Delta^\sigma_{f,i}$ is homeomorphic to $\mathbb{S}^{n-1}\times \mathbb{S}^1$. 

   Without loss of generality suppose that the edge $e_\sigma$ corresponds to the union of stable manifolds of points belonging to the orbit $\mathcal{O}_{\sigma}$ of the point $\sigma\in \Omega^{1}_f$. Denote by $m_\sigma$ the  period of the point $\sigma$.

The set  $cl~W^s_{\mathcal{O}_\sigma}$  cuts the sphere  $S^{n-1}$ in $m_\sigma+1$  connected components that belongs  to  two    $f$-invariant sets $D_{\sigma,-}$ and $D_{\sigma,+}=S^n\setminus cl~D_-$.  Suppose that a domain  correspondent to the vertex $v_\sigma$ belongs to $D_-$. Then $\Delta^{\sigma}_{f,i}=cl\,(p_f(D_{\sigma,-}\setminus (A_f\cup R_f)))$   and the set   $D_{\sigma,-}$ is a union of open balls $D_0, D_1=f(D_0),\dots, D_{m_\sigma-1}=f^{m_\sigma-1}(D_0)$ bounded by $cl\,W^s_{\mathcal{O}_\sigma}$.

From other hand, the set $D_{\sigma,-}$ is a union of stable manifolds of periodic points belonging to the set  $A_{f}\cap D_{\sigma,-}$. Hence, due to Proposition~\ref{ind} the set  $\widehat{D}_{\sigma,-}=p_f(D_{\sigma,-}\setminus (A_f\cup R_f))$ is homeomorphic to $\mathbb{S}^{n-1}\times \mathbb{S}^1$.
	
	Denote by  $l_{\sigma,-}$, $l_{\sigma,+}$   the one-dimensional separatrices   of the point $\sigma$ such that   $l_{\sigma,-}\subset D_{\sigma,-}$, $l_{\sigma,+}\subset D_{\sigma,+}$. Let $N_\sigma$  be the linearizing neighborhood of the point  $\sigma$,   $N^s_\sigma=N_\sigma\setminus W^s_\sigma$, $N^u_\sigma=N_\sigma\setminus W^u_\sigma$,  $\widehat{N}^u_\sigma=N^u_\sigma/_f$, $\widehat{N}^s_\sigma=N^s_\sigma/_f$.

It follows from Propositions~\ref{st-nbh},\ref{grgupo-ado} that the set   $cl~\widehat{D}_{\sigma,-}\cap \widehat{N}^s_\sigma$ is  homeomorphic to the direct product  $\mathbb{S}^{n-2}\times \mathbb{S}^1\times [0,1]$. Then the set     $cl~\widehat{{D}}_{\sigma,-}$  is homeomorphic to the set  $\widehat{{D}}_{\sigma,-}\setminus int~\widehat{N}^s_\sigma$.     
	
	
			The set  $\hat{l}_{\sigma,-}=l_{\sigma,-}/_{f}$ is the  essential knot in $\widehat{D}_{\sigma,-}$ and the set $\widehat{N}^u_{\sigma,-}=(N^u_\sigma\cap D_{\sigma,+})/_f$ is its tubular neighborhood. Due to Corollary~\ref{blue}  the set  $\widehat{D}_+\setminus~int\,\widehat{N}^u_{\sigma,+}$ is homeomorphic to $\mathbb{B}^{n-1}\times \mathbb{S}^1$. Since $\widehat{D}_+\setminus~int\,\widehat{N}^u_{\sigma,+}=\widehat{{D}}_{\sigma,-}\setminus int~\widehat{N}^s_\sigma$ then $\widehat{{D}}_{\sigma,-}\setminus int~\widehat{N}^s_\sigma$ and   $cl~\widehat{{D}}_{\sigma,-}$  are also homeomorphic to $\mathbb{B}^{n-1}\times \mathbb{S}^1$.      
			
			 Now suppose that the  graph $\widehat{\Gamma}_f$ has a loop. Then arguments above prove Lemma for $i\in \{0,\dots,s-1\}$ and the  set   $\Delta_{f,s}=cl\,(\widehat{V}_f\setminus int\,\widehat{N}_{\sigma_*})$, where $N_{\sigma_*}$ is the linearizing neighborhood of the point $\sigma_*$ of negative orientation type and $\widehat{N}_{\sigma_*}=N_{\sigma_*}/f$. Prove that $\Delta_{f,s}$ is homeomorphic to $\mathbb{S}^{n-1}\times \mathbb{S}^1$.

Without loss of generality suppose that   $\sigma_*\in \Omega^{1}_f$.  The set $cl~W^s_{{\sigma_*}}$ cuts the sphere  $S^n$ into two connected components  $D, f(D)$ of period  2. Hence the point ${\sigma_*}$ cuts the set $A_f$ into two symmetrical parts of period 2. Denote by $A_+, A_-$ connected components of the set $A_f\setminus W^u_{\sigma_*}$ lying in the sets  $D, f(D)$ correspondingly and by $l^u_+\subset D, l^u_-\subset f(D)$ unstable separatrices  of the point ${\sigma_*}$.  Put $D_+=D\setminus (A_+\cup R_f)$, $D_-=f(D)\setminus (A_-\cup R_f)$.

It follows from Proposition~\ref{ind}   that  orbit spaces $\widehat{D}_+=D_+/{f^2}$,  $\widehat{D}_-=D_-/{f^2}$ are   homeomorphic  to  $\mathbb{S}^{n-1}\times \mathbb{S}^1$. Since one-dimensional  separatrices of the point ${\sigma_*}$ are fixed with respect $f^2$  then their projections in $\widehat{D}_+, \widehat{D}_-$ are essential knots.    Due to Propositions~\ref{st-nbh},\ref{nb-equ}  the sets    $\widehat{D}_+\setminus~int~\widehat{N}^u_{{\sigma_*}}$, $\widehat{D}_-\setminus~int~\widehat{N}^u_{{\sigma_*}}$ are homeomorphic to $\mathbb{B}^{n-1}\times \mathbb{S}^1$. Then the set $D_+\setminus N^u_{\sigma_*}$ is homeomorphic to $\mathbb{B}^{n-1}\times \mathbb{R}$ and it is possible to find a fundamental domain $B_+\subset D_+\setminus N^u_{\sigma_*}$ of action of $f^2$ on $D_+\setminus N^u_{\sigma_*}$ homeomorphic to $\mathbb{B}^{n-1}\times [0,1]$. Since sets $D_+\cup D_-$ and  $N^u_{\sigma_*}$ are $f$-invariant and 2-periodic, then the domain $B_+$ is also the fundamental domain of action of $f$ on $(D_+\cup D_-) \setminus N^u_{\sigma_*}$  and the factor-space $((D_+\cup D_-) \setminus N^u_{\sigma_*})/_f$ is homeomorphic to $\mathbb{B}^{n-1}\times \mathbb{S}^1$.

Since $V_f\setminus int\,N^s_{\sigma_*}=(D_+\cup D_-)\setminus int\,N^u_{\sigma_*}$ and  $cl\, (\widehat{V}_f\setminus  \widehat{N}^s_{\sigma_*})=\widehat{V}_f\setminus  int\,\widehat{N}^s_{\sigma_*}$ then $cl\, (\widehat{V}_f\setminus  \widehat{N}^s_{\sigma_*})=((D_+\cup D_-) \setminus N^u_{\sigma_*})/_f$, so, $cl\, (\widehat{V}_f\setminus  \widehat{N}^s_{\sigma_*})$ is  homeomorphic to $\mathbb{B}^{n-1}\times \mathbb{S}^1$.

 	\section{A linear-time algorithm for distinguishing edge-colored trees, equipped with automorphisms}\label{dima0}

 In this section we state   the existence of a linear-time algorithm for distinguishing  colored graphs of the cascades from the class $G(S^n)$ that proves Theorem~\ref{dima}.

At first, we will recall some basic definitions and present notation to be used.

The following proposition is  well known in graph theory. 

\begin{st}\label{tr} A connected graph  with $k+1$ vertices  is a tree if and only if  it has $k$ edges. 
\end{st}
\begin{demo}
First prove that the tree  has at least one hanging vertex.   Take an arbitrary vertex $v$ of the tree and move by  an arbitrary edge $e$ incident to $v$ to another vertex $w$. If $w$ is incident only to $e$ then stay at $w$. If there is an edge $l$ incident to $w$ and different from $e$ then go by $l$ to a third vertex and continue the process. Since the tree does not contain cycles we never reach  the same vertex twice. Since the number of the vertices is finite once   the traveling will  stop at a hanging vertex. 

Remove from the tree a hanging vertex and the edge  incident to this vertex.  The obtaining graph is also the tree, so in has a hanging vertex. Remove this vertex and the edge incident to it.   After $k$ iteration we obtain a graph consisting of exactly one vertex and no edges. Since during each iteration  exactly one edge was deleted than originally it was $k$ edges. So, any tree with $k+1$ vertices has exactly $k$ edges.   

If connected graph with $k+1$ vertices and $k$ edges is not a tree then in contains cycles.  Then we can remove  from each cycle an  arbitrary edge to get a tree with  the same number of vertices but with fewer number of edges, that contradicts to the previous  paragraph. 
\end{demo}

Let $k$ be a fixed natural number. Let $T$ be a tree and~\mbox{$c:~E(T)\longrightarrow [k]$}, where $[k]=\{1,2,\ldots,k\}$, be some mapping, called an \emph{edge $k$-coloring of $T$}.
Elements of the set $[k]$ are called \emph{colors}. The pair $(T,c)$ is said to be an \emph{edge-$k$-colored tree}. Let $P$ be some \emph{automorphism} of $(T,c)$, i.e. an
automorphism of~$T$, such that, for any edge $uv\in E(T)$, we have $c(uv)=c(P(u)P(v))$. We will refer to the triple $(T,c,P)$ as an \emph{equipped tree}. Two equipped trees~$(T_1,c_1,P_1)$ and $(T_2,c_2,P_2)$ will said to be \emph{isomorphic}, if there is an isomorphism $\xi$ between $T_1$ and $T_2$, keeping edge colors, i.e. $\forall uv\in E(T_1)~[c(uv)=c(\xi(u)\xi(v))]$), and
conjugating $P_1$ and $P_2$, i.e. $\xi P_1=P_2\xi$.

Meeting the equality $|V(T_1)|=|V(T_2)|$ is a necessary condition for isomorphism of $(T_1,c_1,P_1)$ and $(T_2,c_2,P_2)$. Further, we will consider that both trees $T_1$ and~$T_2$ have
$k$ vertices. We will also assume that $P_1$ and $P_2$ are given by tables, in which, for any vertex $v\in V(T_i)$, a vertex $P_i(v)$, $i=\overline{1,2}$ is also given. We consider that
$c_1$ and $c_2$ are given by tables. Finally, we suppose that $T_1$ and $T_2$ are stored by adjacency lists, i.e. all neighbors are listed, for any vertex of the trees. It will be shown that
the isomorphism problem for two~$k$-vertex equipped trees can be solved in $O(k)$ time.

Let $(T,c,P)$ be some equipped tree. Recall that in Section~\ref{start1} the  notion of the rank of the vertex was given in the following way.   We associate with $T$ a sequence~$T_0,T_1,$ $\ldots,T_s$ of trees, such that $T_0=T$, $T_s$ contains one or two vertices and, for any $i\in [s]$, a tree $T_i$ is obtained
from $T_{i-1}$ by deletion of all its \emph{leaves}, i.e. degree one vertices. All the vertices of $T_s$ are called \emph{central} vertices of the tree $T$ and if $T_s$ has an edge then it is called {\it the central} edge of the tree $T$. The tree $T$ will be \emph{central}, if it has exactly one central
vertex, and \emph{bicentral}, otherwise. \emph{The rank of a vertex $x\in V(T)$}, denoted by $rank(x)$, is the number $\max\{i|~x\in V(T_i)\}$.

By the equipped tree $(T,c,P)$, we will construct a weighted edge-$k$-colored tree $(\hat{T},c,w)$. Vertices of $\hat{T}$ are all the orbits of $P$, two vertices of $\hat{T}$ are connected by an edge if
they are neighbour in $T$. As the weight of a vertex of the tree $\hat{T}$, we take the number of elements in the corresponding orbit of~$P$. As the color of an edge
$O'O''$ of $\hat{T}$, we take the color of any edge of~$T$, simultaneously incident to a vertex from $O'$ and to a vertex from $O''$. It is easy to see that
by the equipped tree the tree $(\hat{T},c,w)$ can be uniquely constructed, see the third part of Lemma~\ref{dima-l1}.

It is easy to see that if $T$ is bicentral and $P$ maps every vertex of $T$ into itself, then the set of central vertices of $T$ coincides
with the set of central vertices of $\hat{T}$. Thereby, knowing $\hat{T}$, it is possible to uniquely restore the set of central vertices of $T$.
Hence, by the second and third parts of Lemma~\ref{dima-l1}, by $(\hat{T},c,w)$ the equipped tree $(T,c,P)$ can be uniquely restored.

Let us construct a simple graph $G$ by $(\hat{T},c,w)$ in the following way. The operation of \emph{$s$-subdivision} of some edge $xy$
of a graph is to delete $xy$, add vertices $z_1,z_2,\ldots,z_s$ and edges $xz_1,z_1z_2,z_2z_3,\ldots,z_{s-1}z_s,z_sy$. The operation of \emph{joining a $s$-cycle to a vertex} $v$ of
some graph is to add vertices $u_1,\ldots,u_{s-1}$, and edges $vu_1,u_1u_2,\ldots,u_{s-2}u_{s-1},u_{s-1}v$ to the graph. For any vertex $v$, we join a cycle of length $w(v)+2$, where
$w(v)$ is the weight of $v$. For any edge~$e\in E(\hat{T})$, we apply its $c(e)$-subdivision, where $c(e)$ is the color of $e$. The resultant
graph is $G$. In the following figure, a weighted colored tree and the corresponding graph are depicted, where we apply
1-subdivision to green edges, 2-subdivision to the blue ones, 3-subdivision to red edges.

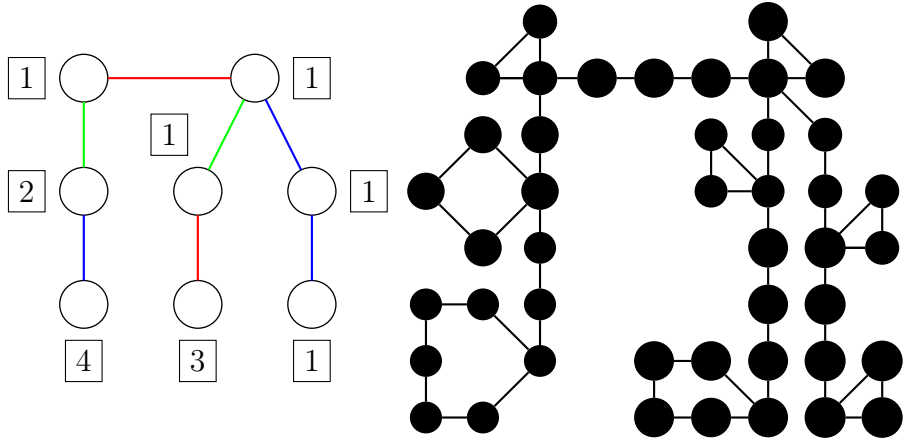
\begin{figure}
\begin{center}
\begin{tikzpicture}[scale=.75]
\Vertex[NoLabel,x=6,y=5]{A}
\Vertex[NoLabel,x=6,y=3]{C}
\Vertex[NoLabel,x=6,y=1]{F}

\Vertex[NoLabel,x=9,y=5]{B}
\Vertex[NoLabel,x=8,y=3]{D}
\Vertex[NoLabel,x=10,y=3]{E}

\Vertex[NoLabel,x=8,y=1]{G}
\Vertex[NoLabel,x=10,y=1]{H}

\Edges[color=red](A,B)
\Edges[color=red](D,G)
\Edges[color=green](A,C)
\Edges[color=green](B,D)
\Edges[color=blue](C,F)
\Edges[color=blue](B,E)
\Edges[color=blue](E,H)

\tikzset{VertexStyle/.style = {
shape = rectangle,
draw}}

\Vertex[x=5,y=5]{1}
\Vertex[x=10,y=5]{1}

\Vertex[x=5,y=3]{2}
\Vertex[x=7.5,y=4]{1}
\Vertex[x=11,y=3]{1}

\Vertex[x=6,y=0]{4}
\Vertex[x=8,y=0]{3}
\Vertex[x=10,y=0]{1}

\tikzset{VertexStyle/.style = {
shape = circle,
fill=black,
inner sep = 0pt,
outer sep = 0pt,
minimum size = 10pt,}}

\Vertex[x=14,y=5]{$a_1$}
\Vertex[x=13,y=5]{$a_2$}
\Vertex[x=14,y=6]{$a_3$}
\Edges($a_1$,$a_2$,$a_1$,$a_3$,$a_2$)

\Vertex[x=14,y=4]{$b_1$}
\Vertex[x=14,y=3]{$b_2$}
\Vertex[x=13,y=4]{$b_3$}
\Vertex[x=13,y=2]{$b_4$}
\Vertex[x=12,y=3]{$b_5$}
\Edges($a_1$,$b_1$,$b_2$,$b_3$,$b_5$,$b_4$,$b_2$)

\Vertex[x=14,y=2]{$c_1$}
\Vertex[x=14,y=1]{$c_2$}
\Vertex[x=14,y=0]{$c_3$}
\Vertex[x=13,y=1]{$c_4$}
\Vertex[x=12,y=1]{$c_5$}
\Vertex[x=12,y=0]{$c_6$}
\Vertex[x=12,y=-1]{$c_7$}
\Vertex[x=13,y=-1]{$c_8$}

\Edges($b_2$,$c_1$,$c_2$,$c_3$,$c_4$,$c_5$,$c_6$,$c_7$,$c_8$,$c_3$)

\Vertex[x=15,y=5]{$d_1$}
\Vertex[x=16,y=5]{$d_2$}
\Vertex[x=17,y=5]{$d_3$}
\Vertex[x=18,y=5]{$d_4$}
\Vertex[x=18,y=6]{$d_5$}
\Vertex[x=19,y=5]{$d_6$}

\Edges($a_1$,$d_1$,$d_2$,$d_3$,$d_4$,$d_5$,$d_6$,$d_4$)

\Vertex[x=18,y=4]{$e_1$}
\Vertex[x=18,y=3]{$e_2$}
\Vertex[x=17,y=4]{$e_3$}
\Vertex[x=17,y=3]{$e_4$}

\Edges($d_4$,$e_1$,$e_2$,$e_3$,$e_4$,$e_2$)

\Vertex[x=18,y=2]{$f_1$}
\Vertex[x=18,y=1]{$f_2$}
\Vertex[x=18,y=0]{$f_3$}
\Vertex[x=18,y=-1]{$f_4$}
\Vertex[x=17,y=0]{$f_5$}
\Vertex[x=16,y=0]{$f_6$}
\Vertex[x=16,y=-1]{$f_7$}
\Vertex[x=17,y=-1]{$f_8$}

\Edges($e_2$,$f_1$,$f_2$,$f_3$,$f_4$,$f_5$,$f_6$,$f_7$,$f_8$,$f_4$)

\Vertex[x=19,y=4]{$g_1$}
\Vertex[x=19,y=3]{$g_2$}
\Vertex[x=19,y=2]{$g_3$}
\Vertex[x=20,y=2]{$g_4$}
\Vertex[x=20,y=3]{$g_5$}

\Edges($d_4$,$g_1$,$g_2$,$g_3$,$g_4$,$g_5$,$g_3$)

\Vertex[x=19,y=2]{$h_1$}
\Vertex[x=19,y=1]{$h_2$}
\Vertex[x=19,y=0]{$h_3$}
\Vertex[x=19,y=-1]{$h_4$}
\Vertex[x=20,y=0]{$h_5$}
\Vertex[x=20,y=-1]{$h_6$}

\Edges($g_4$,$h_1$,$h_2$,$h_3$,$h_4$,$h_5$,$h_6$,$h_4$)
\end{tikzpicture}
\caption{A colored weighted tree and its corresponding graph}
\label{5}
\end{center}
\end{figure}

It is not hard to see that $G$ can be uniquely obtained by $(\hat{T},c,w)$. Conversely,
by $G$ the triple $(\hat{T},c,w)$ can also be restored in unique way. Indeed,
vertices of $G$ of degree more than two correspond to vertices of $\hat{T}$, lengths of cycles are equal to
weights of the corresponding vertices minus two. Lengths of paths of $G$, whose end vertices are of degree
at least three and all internal vertices have degree two, define colors of edges of $\hat{T}$.

Clearly, the graph $G$ is \emph{planar}, i.e. a simple graph, which can be drawn on the plane, such that its vertices are points of the plane and edges are Jordan curves, not intersecting
in internal points. There is known an algorithm with the complexity $O(n)$ for distinguishing $n$-vertex planar graph [1]. Our linear-time algorithm to recognize isomorphism of the equipped trees
$(T,c,P)$ and $(T',c',P')$ is based on this fact.

\begin{Lemm}  One can consider that each of the trees $T$ and $T'$ is bicentral and each of the  automorphisms $P, P'$ maps each of the central  into itself.
\end{Lemm}
\begin{demo} Since $T$ and $T'$ are given by adjacency lists, then the sets of their central vertices
can be computed in $O(n)$ time. Hence, one can consider that $T$ and $T'$ are either simultaneously
central or simultaneously bicentral, otherwise they would not be isomorphic. Assume that
$T$ and~$T'$ are simultaneously bicentral. Denote by $v_1$ and $v_2$ the central vertices of~$T$, and denote by $u_1$ and $u_2$
the central vertices of $T'$. One can consider that $c(v_1v_2)=c'(u_1u_2)$, otherwise the equipped trees are not isomorphic.
Similarly, we may assume that $P(v_1)=v_1, P(v_2)=v_2, P'(u_1)=u_1, P'(u_2)=u_2$ or
$$P(v_1)=v_2, P(v_2)=v_1,P'(u_1)=u_2, P'(u_2)=u_1.$$ In the first situation, we have a case from the statement. In the second one, we define mappings
$\tilde{P}$ and $\tilde{P}'$ in the following way: $$\tilde{P}(v_1)=v_1, \tilde{P}(v_2)=v_2, \forall v\in V(T)\setminus\{v_1,v_2\}~[\tilde{P}(v)=P(v)],$$ 
$$\tilde{P}'(u_1)=u_1,\tilde{P}'(u_2)=u_2,\forall u\in V(T')\setminus\{u_1,u_2\}~[\tilde{P}'(u)=P(u)].$$
It is not hard to see that $\tilde{P}$ and $\tilde{P}'$ are isomorphisms of $(T,c)$ and $(T',c')$, correspondingly. The trees
$(T,c,P)$ and $(T',c',P')$ are isomorphic if and only if $(T,c,\tilde{P})$ and $(T',c',\tilde{P}')$  are isomorphic. Let us consider the case, when
both $T$ and $T$ have exactly one central vertex. Clearly each of  the isomorphisms $P,P'$ maps the central vertex  into
itself. We construct colored trees $(\tilde{T},\tilde{c})$ and $(\tilde{T}',\tilde{c}')$ and some their automorphisms $\tilde{P}$ and $\tilde{P}'$. To construct $\tilde{P}T$,
we take two copies of $(T,c,P)$, whose central vertices are denoted by $v_1$ and $v_2$. Connect them by an edge, and color it in the first color.
We obtain the tree $(\tilde{T},\tilde{c},\tilde{P})$. The tree $(\tilde{T}',\tilde{c}',\tilde{P}')$ is defined by analogy. Both these trees can be obtained
in linear time on $n$. Clearly that $(\tilde{T},\tilde{c},\tilde{P})$ and $(\tilde{T}',\tilde{c}',\tilde{P}')$ are isomorphic if and only if $(\tilde{T},\tilde{c},\tilde{P})$ and $(\tilde{T}',\tilde{c}',\tilde{P}')$ are isomorphic. This finishes the proof of this lemma.\end{demo}

\medskip
{\bf Proof of Theorem~\ref{dima}}.

We suppose that  trees $\Gamma_f$ and $\Gamma_{f'}$ are  bicentral and  the automorphisms $P_f, P_{f'}$ map each of the central vertices into itself. Using the table representation of $P_f$,
the set of all its orbits can be computed in linear time. To this end, we take an arbitrary vertex $v\in V(\Gamma_f)$, compute~\mbox{$P_f(v),P_f(P_f(v)),\ldots$} until $P_f^{r}(v)\neq v$ and, thereby, we find the orbit of $P_f$, containing $v$. By all the orbits of $P_f$, the adjacency list of $\Gamma_f$, and the mapping~$c_f$,
the triple $(\tilde{\Gamma}_f,\tilde{c}_f,\tilde{P}_f)$ can be computed in linear time. By this triple, the planar graph $G_f$ can be computed in $O(|V(G_f)|)$ time. In any tree, the edges number
is less by one than the number of vertices. Clearly that
$$|V(G_f)|\leq |V(\tilde{\Gamma}_f)|+k\cdot|E(\tilde{\Gamma}_f)|+\sum\limits_{v\in V(\tilde{\Gamma}_f)}(\tilde{w}(v)+1)\leq$$
$$\leq |V(\tilde{\Gamma}_f)|+k\cdot |V(\tilde{\Gamma}_f)|+n+|V(\tilde{\Gamma}_f)|\leq (k+3)\cdot n.$$ The equipped trees $(\Gamma_f,c_f,P_f)$ and $(\Gamma_{f'},c_{f'},P_{f'})$ are
isomorphic if and only if the graphs $G_f$ and $G_{f'}$ are isomorphic. Since $k$ is fixed, then for recognizing isomorphism of $G_f$ and $G_{f'}$ it is possible to use
an algorithm from~\cite{Ho}, having the complexity $O(n)$. Therefore the  theorem~\ref{dima} is true.

\end{document}